\documentclass[aop,MSNbibl,nameyear,dvips]{arximspdf}

\doi{10.1214/11-AOP653}
\volume{40}
\issue{4}
\pubyear{2012}
\firstpage{1483}
\lastpage{1534}

\makeatletter

\setattribute{copyright}{owner}{Public Domain}

\newcommand{\binom}[2]{{#1\choose #2}}

\newtheorem{theorem}{Theorem}[section]
\newtheorem{lemma}[theorem]{Lemma}
\newtheorem{coro}[theorem]{Corollary}

\newproclaim{asn}[theorem]{Assumption}

\newcommand{\R}{\mathbb R}
\newcommand{\N}{\mathbb N}
\newcommand{\E}{\mathbb E}
\renewcommand{\P}{\mathbb{P}}
\newcommand{\LA}{\mathcal{L}}
\newcommand{\FA}{\mathcal{F}}
\newcommand{\DA}{\mathcal{D}}
\newcommand{\BA}{\mathcal{B}}
\newcommand{\eps}{\varepsilon}
\newcommand{\al}{\alpha}
\newcommand{\be}{\beta}
\newcommand{\Ga}{\Gamma}
\newcommand{\De}{\Delta}
\newcommand{\sg}{\sigma}
\newcommand{\io}{\iota}
\newcommand{\Lm}{\Lambda}
\newcommand{\Om}{\Omega}
\newcommand{\vp}{\varphi}

\renewcommand{\L}{\LA^\lambda_\eps(u,T)}
\newcommand{\pa}{\partial}
\newcommand{\mn}{\wedge}
\newcommand{\emp}{\varnothing}
\newcommand{\sst}[2]{\mathop{\sum_{#1,#2 = 1}}_{#1\neq #2}}
\newcommand{\sstl}[2]{\mathop{\sum_{#1,#2 = 1}}_{#1< #2}}
\newcommand{\G}{G^{\lambda,u}_\eps}
\newcommand{\Gr}{G^{\lambda,u}}
\newcommand{\Rt}{\Rightarrow}

\makeatother

\begin{document}
\begin{frontmatter}

\title{Generalized self-intersection local time for a~superprocess over a stochastic flow}
\runtitle{GSILTSSF}

\begin{aug}
\author[A]{\fnms{Aaron} \snm{Heuser}\corref{}\ead[label=e1]{heuseram@cc.nih.gov}}
\runauthor{A. Heuser}
\affiliation{National Institutes of Health}
\address[A]{Epidemiology and Biostatistics\\
National Institutes of Health\\
Mark O. Hatfield Clinical Research Center\\
Rehabilitation Medicine Department\\
Bethesda, Maryland\\
USA\\
\printead{e1}} 
\end{aug}

\received{\smonth{10} \syear{2010}}
\revised{\smonth{2} \syear{2011}}

%
\begin{abstract}
This paper examines the existence of the self-intersection local time
for a superprocess over a stochastic flow in dimensions $d\leq3$,
which through constructive methods, results in a Tanaka-like
representation. The superprocess over a stochastic flow is a
superprocess with dependent spatial motion, and thus Dynkin's proof of
existence, which requires multiplicity of the log-Laplace functional,
no longer applies. Skoulakis and Adler's method of calculating moments
is extended to higher moments, from which existence follows.
\end{abstract}

%
\begin{keyword}[class=AMS]
\kwd[Primary ]{60J68}
\kwd{60G57}
\kwd[; secondary ]{60H15}
\kwd{60J80}.
\end{keyword}
\begin{keyword}
\kwd{Superprocess}
\kwd{stochastic flow}
\kwd{self-intersection}
\kwd{local time}.
\end{keyword}

\end{frontmatter}

\section{Introduction}
Superprocesses (or critical branching particle systems), originally
studied by Watanabe (\citeyear{watanabe}) and Dawson
(\citeyear{dawson2,dawson3}) were first shown by Dynkin
(\citeyear{dynkin}) to have a self-intersection local time (SILT). In
particular, Dynkin was able to show existence of the self-intersection
local time for super Brownian motion in $\R^d$, $d\leq7$, provided the
SILT is defined over a region that is bounded away from the diagonal.
When the region contains any part of the diagonal, through
renormalization, the SILT for super Brownian motion has been shown by
Adler and Lewin (\citeyear{lewin}) to exist in $d\leq3$, and further
renormalization processes have been found to establish existence in
higher dimensions by Rosen (\citeyear{rosen}) and Adler and Lewin
(\citeyear{adler}). In regards to non-Gaussian superprocesses, the
SILT has been shown to exist for certain $\alpha$-stable processes by Adler
and Lewin (\citeyear{adler}), and more recently, encompassing
more $\alpha$ values, by Mytnik and Villa (\citeyear{mytnik}). Of important
note, as the $L^2$-limit of an appropriate approximating process, Adler
and Lewin have shown the existence of a class of renormalized SILTs
(indexed on $\lambda>0$) for the super Brownian motion in dimensions
$d=4$ and $5$ and for the super $\al$-stable processes for
$d\in[2\al,3\al)$. As one removes Dynkin's restriction of bounding away
from the diagonal, a singularity arises from ``local double points''
(i.e., $\mu_s\times\mu_t$ where $t=s$) of the process; cf. Adler and
Lewin (\citeyear{lewin}).\vadjust{\goodbreak} The true self-intersection local time should
not be concerned with such local double points, and thus a heuristic
approach to renormalization is naturally observed in the construction.
It should be noted that though this is the method used in Adler and
Lewin (\citeyear{adler,lewin}), a quite different method for
renormalization was developed by Rosen (\citeyear{rosen}). Both
methods are legitimate renormalizations, and lead to existence in
equivalent dimensions, but for this paper, due to the natural
occurrence of the term involving local double points, the initial of
the two methods will be employed. Moreover, the real beauty of this
constructive proof of existence, as seen in Adler and Lewin
(\citeyear{adler,lewin}), is that the aforementioned approximating
process is ``Tanaka-like'' in form. Thus the limit gives a (quite
simple) ``Tanaka-like'' representation for the renormalized SILT.

Quite often, as in the case of Skoulakis and Adler (\citeyear{skoulakis}),
interaction occurs between particles within the system.
Thus, a major drawback in each of the previous superdiffusions is the
requirement of independent spatial motion. Existence as a weak limit of
a branching particle system, and uniqueness as the solution to a
martingale problem, of the superprocess with dependent spatial motion
(SDSM), as a measure-valued Markov process with state space
$M(\hat\R)$,was shown by Wang (\citeyear{wang}). It was later shown
by Dawson, Li and Wang (\citeyear{dawson}) to exist uniquely as a
process in $M(\R)$, and was then extended by Ren, Song and Wang
(\citeyear{ren}) to $M(\R^d)$. Skoulakis and Adler
(\citeyear{skoulakis}) suggested a different model incorporating
dependent spatial motion by replacing the space--time white noise of
Wang's SDSM with a~Brownian flow of homeomorphisms from $\R^d$ to
$\R^d$, which was referred to as a Superprocess over a Stochastic Flow
(SSF).

As of yet, very little work has been done with regard to the
self-intersection local time for superprocesses with dependent spatial
motion. Of important note is the work of He (\citeyear{he}), in which
the existence of the SILT for a superprocess with dependent spatial
motion, similar to the model of Wang, but discontinuous, is shown to
exist in one dimension as a probabilistic limit. Though this was known
to be true, since the local time of the superprocess with dependent
spatial motion was known to exist in one dimension [cf. Dawson, Li and
Wang (\citeyear{dawson})], He was able to give a similar
``Tanaka-like'' representation for the SILT. This paper will
investigate the existence and further properties of a generalized SILT
for the \mbox{$d$-dimensional} SSF, where the generalization refers to the
shift of the support of the Dirac measure away from the origin, to a
point $u\in\R^d$. Note that if $X_t$ is a Markov process, then $Y_t
\triangleq X_t + u$ is a second, dependent Markov Process. The
generalized SILT at~$u$ can be realized as the intersection local time
of the Markov processes~$X_t$ and~$Y_t$.\looseness=-1

\section{Preliminary definitions}

The SSF is constructed as the weak limit of an $\R^d$ branching
particle system. Much of the work that will follow involves
using
properties of the branching particle system, and thus we will briefly
review this construction. This section follows very closely to the work
of Skoulakis\vadjust{\goodbreak} and Adler (\citeyear{skoulakis}), and the reader is
referenced to this work for further questions. We will let
$\overline{\R}{}^d=\R^d\cup\{\De\}$ denote the one-point (Alexandroff)
compactification of $\R^d$, where $\De$ denotes the ``cemetery.'' We
extend measurable functions $\phi\in\BA(\R^d)$ to $\BA(\overline{\R}
{}^d)$ by setting $\phi(\De)=0$.

Let $\N=\{1,2,\ldots\}$ and set
\[
I\triangleq\bigl\{\al=(\al_0,\al_1,\ldots,\al_N)\dvtx N\geq0, \al_0\in\N,
\al_i\in\{0,1\},
1\leq i\leq N\bigr\},
\]
and for any $\al=(\al_0,\ldots,\al_N)\in I$, let $\vert\al\vert=N$ and
$\al-i=(\al_0,\ldots,\al_{\vert\al\vert-i})$. In addition, we will write
$\al\thicksim_n t$ exactly when $t\in[\frac{\vert\al\vert}{
n},\frac{\vert\al\vert+1}{n})$. Let $M(n)$ be the number of particles
alive at time zero, where the spatial position of each particle is
written as
$(x^n_1,x^n_2,\ldots,x^n_{M(n)})$, and define the initial (atomic)
measure by
\[
\mu_0^{(n)}\triangleq\sum_{i=1}^{M(n)}\delta_{x_i^n}.
\]

For each $n\in\N$, $\{B^{\al,n}\dvtx\al_0\leq M(n),\vert\al\vert=0\}
$ is
defined to be a
family of independent $\R^d$ Brownian motions, stopped at time $t=
n^{-1}$, with $B^{\al,n}_0=x^n_{\al_0}$. A recursive definition then
gives a
tree:
for each $k\in\N$, let $\{B^{\al,n}\dvtx\al_0\leq M(n),\vert\al\vert
= k\}
$ be a
collection of $\R^d$ valued Brownian motions, stopped at time
$t=(\vert\al\vert+1) n^{-1}$, and conditionally independent given the
$\sg$-field generated by $\{B^{\al,n}\dvtx\al_0\leq M(n),\vert\al\vert
< k\}
$ and for
which
\[
B^{\al,n}_t=B^{\al-1,n}_t,\qquad t\leq\vert\al\vert n^{-1}.
\]

In regards to branching, for $n\in\N$ let $\{N^{\al,n}\dvtx\al_0\leq
M(n)\}$ be a
family of i.i.d. copies of $N_n$, where $N_n$ is an $\N$-valued random variable
such that
\[
\P(N_n=k)=\cases{
\frac{1}{2}, &\quad$k=2$,\vspace*{2pt}\cr
\frac{1}{2}, &\quad$k=0$.}
\]

Note that it is implicit in the above that the branching is assumed to be
binary, and that for each $n\in\N$,
\begin{eqnarray*}
\E N_n&=&1,
\\[-2pt]
\E N_n^2-(\E N_n)^2&=&1\vspace*{-2pt}
\end{eqnarray*}
and
\[
\E N_n^q=2^{q-1},\qquad q\in\N.\vspace*{-2pt}
\]
Moreover, it is assumed that the families
$\{B^{\al,n}\dvtx\al_0\leq M(n)\}$ and $\{N^{\al,n}\dvtx\al_0\leq
M(n)\}$ are independent.

The final component is that of the stochastic flow. Let $b\dvtx\R
^d\rightarrow\R^d$ and
$c\dvtx\R^d\rightarrow M(d,m)$, where $M(d,m)$ is the space of $d\times m$
matrices, $m\in
\N$, satisfying the following:

\begin{longlist}
\item
the global Lipschitz condition
\[
\vert b(x)-b(y)\vert+\vert c(x)-c(y)\vert\leq C\vert x-y\vert
\]
for any $x,y\in\R^d$;\vadjust{\goodbreak}

\item the linear growth condition,
\[
\vert b(x)\vert+\vert c(x)\vert\leq C(1+\vert x\vert)
\]
for any $x\in\R^d$;

\item for all $i=1,2,\ldots,d$, $j=1,2,\ldots,m$, $b_i$ and $c_{ij}$ are
bounded with
bounded and continuous first and second partial derivatives.
\end{longlist}

Assume that $t\mapsto F^n_{s,t}(x)$ is the solution of the stochastic
differential equation
\[
dY_t=c(Y_t)\,dW^n_t,\qquad Y_s=x,
\]
for all $t\geq s$ and $x\in\R^d$, where $W^n$ is a $\R^m$-valued Brownian
motion, independent of the families $\{B^{\al,n}\}$ and $\{N^{\al,n}\}
$. This
defines a unique Brownian flow of homeomorphisms from $\R^d\rightarrow
\R^d$
[Skoulakis and Adler (\citeyear{skoulakis})].

Set $a_n=n^{-1}$ and $k_n=kn^{-1}$. Then the tree of Brownian
motions over the flow is given by the family of processes $Y^{\al,
n}$, defined
in the following way: let $\al\thicksim_n k_n$ for some $k\in\N$.
Over the time interval
$[0,k_n+a_n]$, $Y^{\al,n}$ is defined to be the solution of the $d$-dimensional
stochastic differential equation,
\begin{eqnarray*}
dY_t&=&b(Y_t)\,dB^{\al,n}_t+c(Y_t)\,dW^n_t,\\
Y_0&=&x^n_{\al_0}.
\end{eqnarray*}

Note that existence and strong uniqueness of the aforementioned
solution is
ensured due
to the assumed conditions on $b$ and $c$. Now set $Y^{\al, n}_t=Y^{\al,
n}_{k_n+a_n}$ for $t> k_n+a_n$ and note that due to construction,
\[
Y^{\al,n}_t=Y^{\al-1,n}_t
\]
for $0\leq t\leq k_n$, $k\in\N$.

We now define the stopping times $\tau^{\al,n}$ as follows: for each
$\al\in I$,
let
\[
\tau^{\al,n}=\cases{
0, &\quad if $\al_0>K_n$,\vspace*{2pt}\cr
\displaystyle\min\biggl\{\frac{i+1}{ n}\dvtx0\leq i\leq\vert\al\vert
,N^{\al
|_i,n}=0\biggr\},&\quad if not $\emp$ and $\al_0\leq M(n)$,\vspace
*{2pt}\cr
\displaystyle\frac{1+\vert\al\vert}{n}, &\quad otherwise.}
\]

The stopped tree of processes, with branching, is the family of processes~$X^{\al,n}$ defined by
\[
X^{\al,n}_t = \cases{
Y^{\al,n}_t, &\quad$t<\tau^{\al,n}$,\vspace*{2pt}\cr
\De, &\quad$t\geq\tau^{\al,n}$.}
\]

The measure-valued process for the finite system of particles is
\[
\mu_t^{(n)}(U)=\frac{\#\{\al\thicksim_n t\dvtx X^{\al,n}_t\in U\}}{n}
\]
for $U\in\BA(\R^d)$, where for a topological space $E$, $\BA(E)$
denotes the
$\sg$-field of Borel measurable sets in $E$.\vadjust{\goodbreak}

We define the corresponding filtration $\FA^n$ by
\[
\FA^n_t\triangleq\sg(B^{\al,n},N^{\al,n}\dvtx\vert\al\vert
<k)\vee\sg
(W^n_s\dvtx s\leq
t)\vee\sg(B^{\al,n}_s\dvtx s\leq t,\vert\al\vert=k)
\]
for $t\in[k_n,k_n+a_n)$, $k=0,1,\ldots.$

Let $C^k(E)$ be the space of continuous functions on $E$ having continuous
partial derivatives up to order $k$, and for $\phi\in C^k(R^d)$ let
\[
\pa^k_{i_1i_2\cdots i_k}\phi(x)=\biggl(\frac{\pa^k\phi}{\pa
x_{i_1}\,\pa
x_{i_2}\cdots\pa x_{i_k}}\biggr)(x).
\]
For $\phi\in C^2(\R^d)$ define the second-order operators $L$ and
$\Lm$ by
%
%
\begin{equation}\label{eqgen1particle}
(L\phi)(x)=\frac{1}{2}\sum_{i,j=1}^d a_{ij}(x,x)\,\pa^2_{ij}\phi(x)
\end{equation}
and
\[
(\Lm\phi)(x,y)=\sum_{i,j=1}^d\sg_{ij}(x,y)\,\pa_i\phi(x)\,\pa_j\phi(y),
\]
where
\[
a_{ij}(x,y)=\delta_{ij}b_i(x)b_j(y)+\sg_{ij}(x,y)
\]
and
\[
\sg_{ij}(x,y)=\sum_{\ell=1}^mc_{i\ell}(x)c_{j\ell}(y),
\]
$x,y\in\R^d$, $i,j=1,\ldots,d$.

Furthermore, for each $n\in\N$, $\phi\in C^2(\R^{n\times d})$
define the
second-order operator $L^{n}$ by
%
%
\begin{equation}\label{eqgennparticles}
(L^n\phi)(x)=\frac{1}{2}\sum_{p,q=1}^n\sum_{i,j=1}^d
a^{pq}_{ij}(x)\,\pa_{p_i}\,\pa_{q_j}\phi(x),
\end{equation}
where
\[
a^{pq}_{ij}(x)=\delta_{pq}\delta_{ij}b_i(x_p)b_j(x_q)+\sg_{ij}(x_p,x_q),
\]
$x=(x_1,\ldots,x_n)$, $x_p\in\R^d$, $p=1,\ldots,n$, and
\[
\delta_{ij}=\cases{
1,&\quad$i=j$,\cr
0,&\quad$i\neq j$.}
\]
For any operator $A$ on a Banach space $\BA$, such that
$A\phi=\lim_{t\rightarrow0}t^{-1}\{T_t\phi-\phi\}$ for some
semigroup $T_t$, we will
denote by $\DA(A)\subset\BA$ the domain of $A$. That is,
\[
\DA(A)=\Bigl\{\phi\in
\BA\dvtx\lim_{t\rightarrow0}t^{-1}\{T_t\phi-\phi\} \mbox{ exists}\Bigr\},
\]
where the limit is in the strong sense.
\begin{asn}
\label{asn1} For the remainder of this paper, the assumption
will be made that $L$ is uniformly elliptic.\vspace*{-3pt}
\end{asn}

For each $k\in\N$ we will denote by $C^k_0(\R^d)$ the
subspace of functions in~$C^k(\R^d)$ which vanish at infinity.

For any topological space $E$, let $M_F(E)$ denote the space of finite Borel
measures on $E$, $C_E[0,\infty)$ the space of continuous paths in $E$
and for
any $\ell\in\N$, $C_K^\ell(\R^d)$ the subspace of $C^\ell(\R^d)$
for which the elements
have compact support.

Endow $D_{M_F(\R^d)}[0,\infty)$ with the topology of weak
convergence,\break that
is, $\mu^{(n)}\in D_{M_F(\R^d)}[0,\infty)$ converges to $\mu\in
D_{M_F(\R^d)}[0,\infty)$ provided\break
$\lim_{n\rightarrow\infty}\langle\phi,\mu^{(n)}\rangle=\langle
\phi,\mu\rangle$
for any $\phi\in C_b(\R^d)$,
and let $\Rt$ denote weak convergence. In addition, for any $\mu\in
M_F(E)$ and
$\ell\in\N$, denote by $\mu^\ell$ the product measure
$\mu\times\mu\times\cdot\times\mu\in M_F(\R^{\ell\times d})$.
Under these
assumptions, and Assumption \ref{asn1} upon $L$, we arrive at the
following theorem.\vspace*{-3pt}
\begin{theorem}
\label{thmsp} Let $\mu^{(n)}$ be defined as above with
$\mu_0^{(n)}\Rt\mu_0$, then \mbox{$\mu^{(n)}\Rt\mu$},
where $\mu\in C_{M_F(\R^d)}[0,\infty)$ is the unique solution of the
following
martingale problem:

For all $\phi\in C^2_K(\R^d)$,
%
%
\begin{equation}\label{eqSPDE}
Z_t(\phi)=\langle\phi,\mu_t\rangle-\langle\phi,\mu_0\rangle
-\int_0^t ds \langle L\phi,\mu_s\rangle
\end{equation}
is a continuous square integrable $\{\FA^\mu_t\}$-martingale such that
$Z_0(\phi)=0$ and has quadratic variation process
%
%
\begin{equation}\label{eqbracket}
\langle Z(\phi)\rangle_t=\int_0^t ds (\langle\phi^2,\mu_s\rangle
+\langle\Lm\phi,\mu_s^2\rangle
).\vspace*{-3pt}
\end{equation}
\end{theorem}
\begin{pf}
See Theorem 2.2.1 in Skoulakis and Adler (\citeyear{skoulakis}).\vspace*{-3pt}
\end{pf}
\begin{asn}
\label{asn2} For the remainder of this work, it will be assumed
that $\mu_0\in M_F(\R^d)$ has compact support and satisfies
\[
\mu_0(dx)\leq
m(x)\,dx
\]
for some bounded $m\in L^1(\R^d).$\vspace*{-3pt}
\end{asn}

\section{Some needed lemmata}
Some needed technical lemmata, will be presented, where due to the significantly
large number of calculations required, the
proof is deferred to the \hyperref[app]{Appendix}.

As in most existence proofs for self-intersection local time of a superprocess,
higher moments of the superprocess are required; cf. Adler and Lewin
(\citeyear{adler}),
Dynkin (\citeyear{dynkin}). Through finding the first and second
moments of the branching
process, and passing to the limit as $n\rightarrow\infty$ [Skoulakis
and Adler
(\citeyear{skoulakis})] found the first and second moments for the
SSF. A variation of
this method is employed to find higher moments of the SSF.\vadjust{\goodbreak}
\begin{lemma}
If $L^n$ is defined as (\ref{eqgennparticles}), then $L^n$ is the
generator of the diffusion which describes the joint motion of $n$
particles in
the aforementioned branching particle system.
\end{lemma}
\begin{pf}
See the Appendix of Skoulakis and Adler (\citeyear{skoulakis}).
\end{pf}
\begin{lemma}
\label{lmtf}
For each $n\in\N$, there exists a transition
function $q^n_t$
for the Markov process $Y_t=(Y^1_t,\ldots,Y^n_t)$. Furthermore, $\{
Q^n_t\dvtx t\geq0\}$,
defined by
\[
Q^n_t\phi(x)=\int_{\R^d}\phi(y)q^n_t(x,y)
\]
is a strongly
continuous contraction semigroup on $C_0(\R^d)$.
\end{lemma}
\begin{pf}
Since it is assumed that Assumption \ref{asn1} holds for $L$, it
follows that for
each $n\in\N$, Assumption \ref{asn1} also holds for $L^n$. Theorem
5.11 in
Dynkin (\citeyear{dynkin65})
then completes the proof.
\end{pf}

We denote by $C^\infty(\R^d)$ the space of infinitely differentiable
functions on~$\R^d$, by $C^\infty_K(\R^d)$, the subspace of $C^\infty(\R^d)$ of
which the elements
have compact support, by $\DA'(\R^d)$ the space of distributions on
$C^\infty_K(\R^d)$, and by~$D^\al u$ the $\al$th-weak partial derivative of $u$. Note
that a differentiable function will have a weak derivative that agrees
with the
functions derivative, and thus we will at times use a slight abuse in notation
and write the weak derivative as
$D^\al=\pa^{\al_1}_1\,\pa^{\al_2}_2\cdots\pa^{\al_d}_d$.

We denote by $S_d$ the Schwartz space of rapidly decreasing functions
on~$\R^d$,
and the dual to $S_d$, the space of tempered distributions on $R^d$, by
$S'_d$. For any two functions $\phi\dvtx E_1\rightarrow\R$, $\psi
\dvtx E_2\rightarrow\R$ denote by
$\phi\otimes\psi$ the concatenation of $\phi$ and $\psi$. That is,
$\phi\otimes\psi\dvtx E_1\times E_2\rightarrow\R$ is the map defined by
$(x_1,x_2)\mapsto\phi(x_1)\psi(x_2)$.
\begin{lemma}
\label{lmsinglevariableapprox}Let $\phi\in S_{\ell\times d}$,
then there
exists $\{\phi_n\dvtx n\in\N\}$ such that:
\begin{longlist}
\item $\phi_n=\sum_{k=1}^n\phi^1_k\otimes\phi^2_k\otimes
\cdots\otimes
\phi^\ell_k$,
for some $\phi^1_k,\ldots,\phi^\ell_k\in C_K^\infty(\R^d)$;

\item $\phi_n$ converges to $\phi$ in $S_{\ell\times d}$ as
$n\rightarrow\infty$.
\end{longlist}
\end{lemma}
\begin{pf}
Taylor's theorem implies the above holds for any $\phi\in C^\infty
_K(\R^d)$;
cf. Rudin (\citeyear{babyrudin,greenrudin}). From Theorem 7.10 of
Rudin (\citeyear{rudin})
there exist
$\{\phi_n\dvtx n\in\N\}\subset C^\infty_K(\R^{\ell\times d})$ such
that $\phi_n$
converges to $\phi$ in $S_d$, and the result thus follows.
\end{pf}

Given $\phi\in\BA(\R^{(n+1)\times d})$, $n\in\N$, define the projection
$\pi_1$ by
\[
(\pi_1
Q^n_t\phi)(x_1,\ldots,x_{n+1})=Q^n_t\phi_{x_1}(x_{2},\ldots,x_{n+1}),
\]
where
$\phi_{x}(y_1,\ldots,y_n) =\phi(x,y_1,\ldots,y_n)$.\vadjust{\goodbreak}

Given\vspace*{2pt} $m\in\N$, $i=1,2,\ldots,m-1$, $j=1,2,\ldots,m$, $i\neq j$ and any
function $\phi\dvtx\R^{m\times d}\rightarrow\R$, define
$(\Phi_{ij}\phi)\dvtx\R^{(m-1)\times d}\rightarrow\R$ by
\begin{eqnarray*}
&&
(\Phi_{ij}\phi)(x_1,\ldots,x_{m-1})\\
&&\qquad=
\cases{\phi(x_1,\ldots,x_i,\ldots,x_{j-1},x_i,x_j,\ldots,x_{m-1}),
&\quad
$i<j, j\neq m$,\vspace*{2pt}\cr
\phi(x_1,\ldots,x_i,\ldots,x_{m-1},x_i), &\quad$i<j, j=m$,\vspace
*{2pt}\cr
\phi(x_1,\ldots,x_{j-1},x_i,x_j,\ldots,x_i,\ldots,x_{m-1}), &\quad$i>j,
j\neq1$,\vspace*{2pt}\cr
\phi(x_i,x_1,\ldots,x_i,\ldots,x_{m-1}), &\quad$i>j, j=1$.}
\end{eqnarray*}

Furthermore, for $\mathbf{x}_m=(x_1,x_2,\ldots,x_m)$, $x_p\in\{ij,0\}$,
$p=1,2,\ldots,m$, let
\[
\zeta(x_p)=\cases{
\Phi_{ij}, &\quad$x_p=ij$,\vspace*{2pt}\cr
\pi_1, &\quad$x_p=0$,}
\]
and
\[
\ell(x_p)=\cases{
\ell(x_{p-1})+1, &\quad$x_{p-1}=ij$,\vspace*{2pt}\cr
\ell(x_{p-1})-1, &\quad$x_{p-1}=0$,}
\]
$\ell(x_0)=\ell\in\N$. Let
$\mathbf{s}_m=(s_1,s_2,\ldots,s_m)$, $s_p\in[0,\infty)$, $p=1,2,\ldots,m$,
and denote
%
%
\begin{equation}\label{fnGamma}
\Ga^{\mathbf{x}_{m - 1}}_{\ell;\mathbf{s}_{m}} \triangleq
Q^{\ell}_{s_1}\zeta(x_1)Q^{\ell(x_1)}_{s_2-s_1}\zeta(x_2)Q^{\ell
(x_2)}_{s_3-s_2}
\cdots
\zeta(x_{m-2})Q^{\ell(x_{m-2})}_{s_{m - 1} - s_{m - 2}}\pi
_1Q^{\ell(x_{m-1})}_
{ s_m - s_ { m - 1 } }.\hspace*{-28pt}
\end{equation}

The next lemma comes from Skoulakis and Adler (\citeyear{skoulakis}),
though it should be noted that in the aforementioned
paper the result is shown for near critical branching (as opposed to
critical branching in this paper). This is of little concern though, as
a modification of the original proof (making for a much simpler proof)
gives the critical branching case.
\begin{lemma}
\label{lmmoment12sp} Let $\phi,\phi_1,\phi_2\in C^2_K(\R
^d)$ and $t>0$,
then
\[
\mbox{\textup{(i)}}\quad
\E\mu_t(\phi)=\langle Q_t\phi, \mu_0\rangle
\]
and
\begin{eqnarray*}
\mbox{\textup{(ii)}}\quad
\E\mu_{t_1}(\phi_1)\mu_{t_2}(\phi_2)&
=&\bigl\langle Q^2_{t_1}\bigl(\pi_1Q_{t_2-t_1}(\phi_1\otimes\phi_2)\bigr),\mu
_0^2\bigr\rangle\\
&&{}+\int_0^{t_1} ds \bigl\langle Q_s\Phi_{12}Q^2_{t_1-s}\bigl(\pi_1Q_{t_2-t_1}
(\phi_1\otimes\phi_2)\bigr),\mu_0\bigr\rangle
\end{eqnarray*}
with the convention that $Q^n_0\phi=\phi$, $n\in\N$.
\end{lemma}
\begin{pf}
See Skoulakis and Adler (\citeyear{skoulakis}), Proposition 3.2.1
[with $(1+\gamma_n/n)$ and $e^{-\lambda\gamma r_n}$ both replaced by $1$].
\end{pf}

Before our moment calculations, some needed definitions and lemmata
will be
presented. In what follows $(S,d)$ will refer to a metric space, in
which it is
assumed $S$ is separable, and $\rho$ will denote the Prohorov metric on
$M_F(S)$.
\begin{lemma}
\label{lmEKconv2}
If $\{\mu^{(n)}\dvtx n\geq0\}\subset M_F(\R^d)$ satisfies
$\mu^{(n)}\Rt\mu\in M_F(\R^d)$
then
\[
\bigl(\mu^{(n)}\bigr)^\ell\Rt\mu^\ell
\]
for all $\ell\in\N$.
\end{lemma}
\begin{pf}
Define
\[
M=\Biggl\{\phi=\bigotimes_{k=1}^\ell\phi_k\dvtx\ell\geq1,\phi_k\in
C_K(\R^d)\cup\{1\}, k=1,2,\ldots,\ell\Biggr\}.
\]
From Ethier and Kurtz (\citeyear{ethier}) Chapter 3, Proposition
4.4, for any $\nu,\nu^{(n)}\in M_F(\R^d)$, $n=1,2,\ldots,$ such that
$\lim_{n\rightarrow\infty}\langle\phi,\nu^{(n)}\rangle=\langle
\phi,\nu\rangle$ for all $\phi\in
C_K(\R^d)$, it follows that $\nu^{(n)}\Rt\nu$. For any $\ell\in\N
$, since
$\mu^{(n)}\Rt\mu$,
$\lim_{n\rightarrow\infty}\langle\phi,(\mu^{(n)})^\ell\rangle
=\langle\phi,\mu^\ell\rangle$, for any
$\phi=\bigotimes_{k=1}^\ell\phi_k $ with $\phi_k\in C_K(\R^d)$ or
$\phi_k\in\{1\}$, $k=1,\ldots,\ell$. Thus, for any $\phi=\bigotimes
_{k=1}^\ell\phi_k
\in M$, $\lim_{n\rightarrow\infty}\langle\phi,(\mu^{(n)})^\ell
\rangle=\langle\phi,\mu^\ell\rangle$, which
implies, by Ethier and Kurtz (\citeyear{ethier}), Chapter 3,
Proposition 4.6,
$(\mu^{(n)})^\ell\Rt\mu^\ell$.
\end{pf}

We denote by $D_S[0,\infty)$ the Skorohod space on $S$, that is, the space
of all c\`{a}dl\`{a}g mappings from $[0,\infty)$ to $S$. Note that
under the
assumption
that $S$ is separable, $D_S[0,\infty)$ with the
metric defined by Ethier and Kurtz (\citeyear{ethier}), Chapter 3,
(5.2), is a separable metric
space. Moreover, if $(S,d)$ is complete, $D_S[0,\infty)$ is complete; cf.
Ethier and Kurtz (\citeyear{ethier}), Theorem 5.6, Chapter 3. For
$\phi\in C_b(S)$ define the metric
$\Vert\phi\Vert_{bL}=\Vert\phi\Vert_\infty\vee\sup_{x\neq
y}\frac{\vert\phi(x)-\phi(y)\vert}{d(x,y)}$. The next two lemmata are
essential in
the moment proofs for the superprocess.
\begin{lemma}
\label{lmEKconv3} For $k,\ell\in\N$, let $\psi\dvtx\R
^k_+\times\R^\ell\rightarrow\R$ in
$C_b(\R^k_+\times\R^\ell)$ satisfy $\sup_{s\in\R^k_+}\Vert\psi
(s,\cdot)\Vert_{bL}<\infty$, and
let $\mu_0$ be an a.s. finite measure having compact support with
$\mu_0^{(n)}\Rt\mu_0.$ Then
\begin{eqnarray*}
&&\lim_{n\rightarrow\infty}\Biggl|\frac{1}{n^k} \mathop{\sum_{r_1,\ldots,
r_k=0}}_{r_1<\cdots<r_k}^{[nt] -1
} \biggl\langle\psi\biggl(\frac{r}{n}, \cdot\biggr) , \bigl(\mu_0^{(n)}\bigr)^\ell\biggr\rangle\\
&&\qquad\hspace*{2pt}{} - \int_0^t ds_k
\int_0^{s_k} ds_{k - 1} \cdots\int_0^{s_2} ds_1\langle\psi(s,
\cdot) , \mu_0^\ell\rangle
\Biggr|=0,
\end{eqnarray*}
where $r=(r_1,\ldots,r_k)$ and $s=(s_1,\ldots,s_k)$.
\end{lemma}
\begin{pf}
Indeed,
\begin{eqnarray*}
&&
\Biggl|\frac{1}{n^k} \mathop{\sum_{r_1,\ldots,
r_k=0}}_{r_1<\cdots<r_k}^{[nt] -1
} \biggl\langle\psi\biggl(\frac{r}{n}, \cdot\biggr) , \bigl(\mu_0^{(n)}\bigr)^\ell\biggr\rangle
- \int_0^t ds_k
\int_0^{s_k} ds_{k - 1} \cdots\int_0^{s_2} ds_1\langle\psi(s,
\cdot) , \mu_0^\ell\rangle
\Biggr|\\
&&\quad \leq
\frac{1}{n^k}\mathop{\sum_{r_1,\ldots,
r_k=0}}_{r_1<\cdots<r_k}^{[nt] -1
}\biggl\vert\biggl\langle\psi\biggl(\frac{r}{n},\cdot\biggr),\bigl(\mu_0^{(n)}\bigr)^\ell\biggr\rangle
-\biggl\langle\psi\biggl(\frac{r}{n} , \cdot\biggr) , \mu_0^\ell\biggr\rangle\biggr\vert\\
&&\qquad{} +\Biggl\langle\Biggl|\frac{1}{n^k}\mathop{\sum_{r_1,\ldots, r_k=0}}_{r_1<\cdots
<r_k}^{[nt] -1 }\psi\biggl(\frac{r}{n},\cdot\biggr) - \int_0^t ds_k \int
_0^{s_k} ds_{k - 1}\cdots\int_0^{s_2} ds_1\psi(s,\cdot)\Biggr|,\mu
_0^\ell\Biggr\rangle.
\end{eqnarray*}

By assumption $\sup_{s}\Vert\psi(s,\cdot)\Vert_{bL}<\infty$, and
thus from
Ethier and Kurtz (\citeyear{ethier}) the first of the above terms
converges to zero. Since
$\psi$ is continuous and bounded, and $\mu_0^\ell$ is finite with
compact support,
it follows that the second term is also convergent toward zero.
\end{pf}
\begin{lemma}
\label{lmEKconv4}For any $\phi_i\in C^\infty_K(\R^{d})$,
$i=1,2,\ldots,\ell$,
$\ell\in\N$, $0<t<\infty$,
\[
\lim_{n\rightarrow\infty}\E\bigl\langle\phi_1\otimes\phi_2\otimes
\cdots\otimes\phi_\ell, \bigl(\mu_{t}^{(n)}\bigr)^\ell\bigr\rangle=\E
\langle\phi_1\otimes\phi_2\otimes\cdots\otimes\phi_\ell,\mu
_{t}^\ell\rangle.
\]
\end{lemma}
\begin{pf}
Let $\mu^{(n)}=\{\mu_t^{(n)}\dvtx t\geq0\}$ be a branching process as defined
above, let $\mu$ be a weak limit point of $\mu^{(n)}$, and let $\{
n_k\}$ be the
subsequence along which $\mu^{(n_k)}\Rt\mu$. From Ethier and Kurtz
(\citeyear{ethier}), Theorem 3.1, Chapter 3, there is a Skorohod
representation for $\mu,\mu^{(n_k)}$,
$k\in\N$. That is, there exist random variables $X,X_k$, $k\in\N$,
defined on
the same probability space, such that $X\stackrel{d}{=}\mu$,
$X_k\stackrel{d}{=}\mu^{(n_k)}$, $k\in\N$, and $X_k\rightarrow X$
a.s. as
$k\rightarrow\infty$.

For $X\in D_{M_F(\R^d)}[0,\infty)$, define $\P X(\phi_i)^{-1}$ to be the
distribution of $X(\phi_i)\in D_{\R}[0,\infty)$ then, by dominated
convergence
\begin{eqnarray*}
&&\lim_{k\rightarrow\infty}\sup_{\Vert\psi\Vert_{bL}=1}\Biggl\vert
\Biggl\langle\psi,\prod_{i=1}^\ell\P X_k(\phi_i)^{-1}\Biggr\rangle-\Biggl\langle
\psi,\prod_{i=1}^\ell\P X(\phi_i)^{-1}\Biggr\rangle\Biggr\vert\\
&&\qquad=\lim_{k\rightarrow\infty}\sup_{\Vert\psi\Vert_{bL}=1}{|}\E
\psi(X_k(\phi_1),\ldots,
X_k(\phi_\ell))-\E\psi(X(\phi_1),\ldots,X(\phi_\ell)){|}\\
&&\qquad= 0.
\end{eqnarray*}
It then follows from Ethier and Kurtz (\citeyear{ethier}) that
\[
\lim_{k\rightarrow\infty}\rho\Biggl(\prod_{i=1}^\ell\P
X_k(\phi_i)^{-1} , \prod_{i=1}^\ell\P X(\phi_i)^{-1}\Biggr)=0
\]
or equivalently,
\[
(X_k(\phi_1),\ldots,X_k(\phi_\ell))\Rt
(X(\phi_1),\ldots,X(\phi_\ell))
\]
in $D_{\R^\ell}[0,\infty)$. Therefore, from Theorem \ref{thmsp},
\[
\bigl(\mu^{(n_k)}(\phi_1),\ldots,\mu^{(n_k)}(\phi_\ell)\bigr)\Rt(\mu(\phi
_1),\ldots,
\mu(\phi_\ell))
\]
in $D_{\R^\ell}[0,\infty)$. Thus, from Lemma A.3.9 of Skoulakis and
Adler (\citeyear{skoulakis}), for
$i=1,2,\ldots,\mu(\phi_i)$ is continuous. Therefore, the open mapping theorem
[Ethier and Kurtz (\citeyear{ethier}), Chapter 3, Corollary 1.9]
implies that
\[
\bigl(\mu_t^{(n_k)}(\phi_1),\ldots,\mu_t^{(n_k)}(\phi_\ell)\bigr)\Rt(\mu
_t(\phi_1),\ldots,
\mu_t(\phi_\ell))
\]
in $\R^\ell$, which further implies that
\[
\mu_t^{(n_k)}(\phi_1)
\cdot\mu_t^{(n_k)}(\phi_2)\cdots\mu_t^{(n_k)}
(\phi_\ell)\Rt\mu_t(\phi_1)\cdot\mu_t(\phi_2)\cdots
\mu_t(\phi_\ell)
\]
in $\R$. Note that [cf. (3.1) in Skoulakis and Adler (\citeyear
{skoulakis})] for any $t\geq0$,
$\E\mu_t^{(n)}(1)=\mu_0^{(n)}(1)$, and thus $\{\mu_t^{(n)}(1)\dvtx t\geq
0\}$ is an
$\FA^{n}_t$-martingale. It follows from Doob's maximal inequality
[Karatzas and Shreve (\citeyear{shreve}), Theorem 3.8] that for any
$T\geq0$,
\[
\E\sup_{0\leq t\leq
T}\bigl[\mu_t^{(n)}(1)\bigr]^\ell\leq\biggl(\frac{\ell}{\ell
-1}\biggr)^\ell\E\bigl[
\mu_T^{(n)}(1)\bigr]^\ell.
\]
Since $\mu_0^{(n)}\Rt\mu_0$, $\lim_{n\rightarrow\infty} \mu
_0^{(n)}(1)=\mu_0(1)$,
and thus, $\sup_{n\geq1}\mu_0^{(n)}(1)<\infty$. Since $\mu
_t^{(n)}(1)$ is the
total mass process of the branching particle system, and is absent of influence
by the stochastic flow, $[\mu_T^{(n)}(1)]^\ell$ is equivalent in
distribution to a
total mass process with an initial $M(n)^\ell$ particles, which implies
$\E[\mu_T^{(n)}(1)]^\ell=[\mu_0^{(n)}(1)]^\ell$. Thus,
$\sup_{n\geq1}\E\sup_{0\leq t\leq
T}[\mu_t^{(n)}(1)]^\ell<\infty$.
Theorem 25.12 of Billingsley (\citeyear{billingsley}) implies
$\lim_{k\rightarrow\infty}\E\prod_{i=1}^\ell\mu_t^{(n_k)}(\phi
_i)=\E\prod_{i=1}
^\ell\mu_t(\phi_i)$,
and thus,
\[
\lim_{n\rightarrow\infty}\E\bigl\langle\phi_1\otimes\phi_2\otimes
\cdots\otimes\phi_\ell,\bigl(\mu_t^{ (n)}\bigr)^\ell\bigr\rangle=\E\langle
\phi_1\otimes\phi_2\otimes\cdots\otimes\phi_\ell,\mu_t^\ell
\rangle.
\]
\upqed\end{pf}

In Skoulakis and Adler (\citeyear{skoulakis}) the first and second
moment calculations are done via
first finding $\E\langle\phi,\mu_t^{(n)}\rangle$ and
$\E\langle\phi_1\otimes\phi_2,\mu_{t_1}^{(n)}\mu
_{t_2}^{(n)}\rangle$
then passing to the limit as $n\rightarrow\infty$. This works well
when the number of
cases to consider are small, but due to the
rapid growth in cases to consider as the moments increase, the
following method
will vary slightly. The method first
calculates $\E\langle\phi,(\mu_t^{(n)})^3\rangle$ and $\E\langle
\psi,(\mu_t^{(n)})^4\rangle$
for $\phi\in C^\infty_K(\R^{3\times d}),\psi\in C^\infty_K(\R
^{4\times d})$,
$t\geq0$,
then
passes to the limit before utilizing the Markov property to find
$\E\langle\phi_1\otimes\phi_2\otimes\phi_3,\mu_{t_1}\mu
_{t_2}\mu_{t_3}\rangle$ and
$\E\langle\psi_1\otimes\psi_2\otimes\psi_3\otimes\psi_4,\mu
_{t_1}\mu_{t_2}\mu_{t_3 }\mu_{t_4}\rangle$, where $\phi_i,\psi
_j\in C^\infty_K(\R^d)$,
$i=1,2,3$, $j=1,2,3,4$,
and $0<t_1\leq t_2\leq t_3\leq t_4$.

Since the calculations for the third moment are a much simpler case of the
fourth, we
present here only the derivation of the fourth moment. To begin, note that
%
%
\begin{equation}\qquad
\E\bigl\langle\phi,\bigl(\mu_t^{(n)}\bigr)^\ell\bigr\rangle=\frac{1}{n^\ell}\mathop
{\sum
_{\al_k\thicksim_n
t}}_{k=1,\ldots,\ell}\E\phi(Y^{\al_1,n}_t,Y^{\al_2,n}_t,\ldots,Y^{\al
_\ell,n}_t)\E\prod_{
i=1 } ^\ell1_ { \al_i
,n}(t),
\end{equation}
where $1_{\al_i,n}(t)$ is the indicator on the event that the particle
$\al_i$
is alive at time $t$.
Thus, for the fourth moment, if $\al_i\thicksim_n t$, $i=1,2,3,4$ and $N=[t
n]$, we will have the following cases to consider:

\begin{enumerate}[(III)]
\item[(I)] Each particle resides on its own tree.
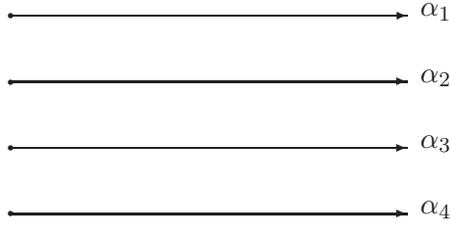
\vspace*{8pt}

\begin{picture}(0,0)
\multiput(0,0)(0,-25){4}{\vector(1,0){150}}
\multiput(0,0)(0,-25){4}{\circle*{2}}
\put(155,0){$\al_1$}\put(155,-25){$\al_2$}\put(155,-50){$\al
_3$}\put(155,
-75)
{$\al_4$}
\end{picture}

\vspace{1.26in}
\item[(II)]
Two particles reside on one tree, the
other two
reside on their own trees. Thus, the two particles on the common tree
share a
common ancestor $\be$ with $\vert\be\vert=r$ and $r\in\{
0,1,\ldots,N-1\}$.
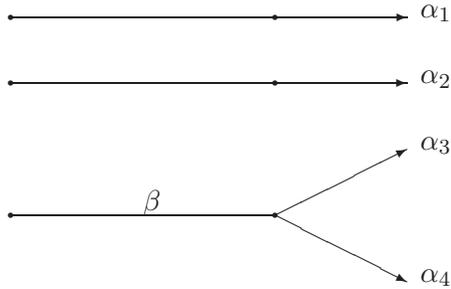
\vspace*{8pt}

\begin{picture}(0,0)
\multiput(0,0)(0,-25){2}{\vector(1,0){150}}
\multiput(0,0)(0,-25){2}{\circle*{2}}
\multiput(100,0)(0,-25){2}{\circle*{2}}
\put(0,-75){\line(1,0){100}}
\multiput(0,-75)(100,0){2}{\circle*{2}}
\put(100,-75){\vector(2,1){50}}
\put(100,-75){\vector(2,-1){50}}
\put(155,0){$\al_1$}\put(155,-25){$\al_2$}\put(155,-50){$\al
_3$}\put(155,-100)
{$\al_4$}\put(50,-73){$\be$}
\end{picture}

\item[(III)]
\vspace{1.51 in}
Two particles reside on one tree, the other two
on a second tree. Thus, the two particles on one tree share a common ancestor
$\be_1$ with $\vert\be_1\vert=r_1$, the two particles on the second tree
have a
common ancestor $\be_2$ with $\vert\be_2\vert=r_2$ and
$r_1,r_2\in\{0,1,\ldots,N-1\}$.
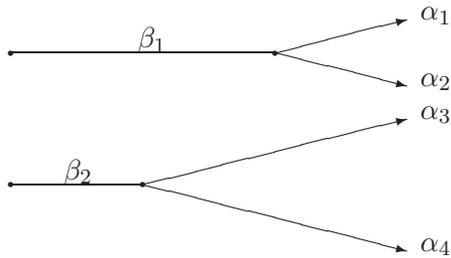
\vspace*{8pt}

\begin{picture}(0,0)
\put(0,-25){\line(1,0){100}}\put(0,-75){\line(1,0){50}}
\put(100,-25){\vector(4,1){50}}\put(100,-25){\vector(4,-1){50}}
\put(50,-75){\vector(4,1){100}}\put(50,-75){\vector(4,-1){100}}
\multiput(0,-25)(100,0){2}{\circle*{2}}\multiput
(0,-75)(50,0){2}{\circle*{2}
}\put(155,-12.5){$\al_1$}\put(155,-37.5){$\al_2$}\put
(155,-50){$\al_3$}\put(155,
-100){$\al_4$}\put(48,-23){$\be_1$}\put(20,-73){$\be_2$}
\end{picture}

\vspace{1.51 in}
\item[(IV)]
Three particles reside on one tree, the
fourth on
its own tree. Thus, two of the three particles share a common ancestor
$\be_2$
with $\vert\be_2\vert=r_2$, and all three share a common ancestor
$\be
_1$ with
$\vert\be_1\vert=r_1$,\vadjust{\goodbreak} such that $r_1\in\{0,1,\ldots,r_2-1\}$ and
$r_2\in\{1,\ldots,N-1\}$.
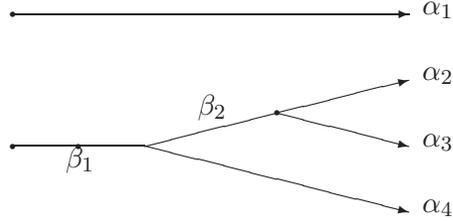
\vspace*{6pt}

\begin{picture}(0,0)
\put(0,0){\vector(1,0){150}}
\put(0,-50){\line(1,0){50}}
\put(50,-50){\vector(4,1){100}}
\put(50,-50){\vector(4,-1){100}}
\put(100,-37.5){\vector(4,-1){50}}
\put(0,0){\circle*{2}}\put(0,-50){\circle*{2}}
\put(25,-50){\circle*{2}}\put(100,-37.5){\circle*{2}}
\put(155,0){$\al_1$}\put(155,-25){$\al_2$}\put(155,-50){$\al
_3$}\put(155,-75)
{$\al_4$}\put(70,-38){$\be_2$}\put(20,-58){$\be_1$}
\end{picture}

\vspace{1.11in}
\item[(V)]
All four particles reside on one tree.
This gives
the following two sub-cases:

\begin{enumerate}[(A)]
\item[(A)]
Two of the particles share a common ancestor $\be_3$ with $\vert
\be_3\vert=r_3$,
the other two share a common ancestor $\be_2$ also with $\vert\be
_2\vert=r_2$, all
four share a common ancestor $\be_1$ with $\vert\be_1\vert=r_1$,
$\be
_2$ and $\be_3$
are both descendants of $\be_1$, and $r_1\in\{0,1,\ldots,(r_2-1)\mn
(r_3-1)\}$ and
$r_2,r_3\in\{1,\ldots,N-1\}$.
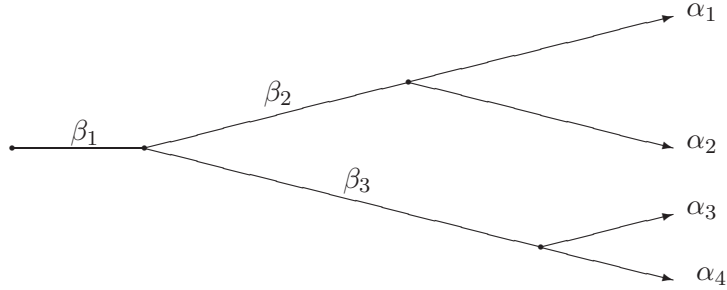
\vspace*{3pt}

\begin{picture}(0,0)
\put(0,-50){\line(1,0){50}}
\put(50,-50){\line(4,1){100}}\put(50,-50){\line(4,-1){150}}
\put(150,-25){\vector(4,1){100}}\put(150,-25){\vector(4,-1){100}}
\put(200,-87.5){\vector(4,1){50}}\put(200,-87.5){\vector(4,-1){50}}
\put(0,-50){\circle*{2}}\put(50,-50){\circle*{2}}\put
(150,-25){\circle*{2}}
\put(200,-87.5){\circle*{2}}
\put(255,0){$\al_1$}\put(255,-50){$\al_2$}\put(255,-75){$\al
_3$}\put(255,-100){
$\al_4$}\put(22,-47){$\be_1$}\put(95,-32){$\be_2$}\put
(125,-65){$\be_3$}
\end{picture}

\vspace{1.51in}
\item[(B)]
Two of the particles share a common ancestor $\be
_3$, another
particle shares a common ancestor $\be_2$ with $\be_3$, all four
particles share
a common ancestor $\be_1$ and $\be_1$ is an ancestor of $\be_2$
which is an
ancestor of $\be_3$, with $\vert\be_1\vert=r_1$, $\vert\be_2\vert=r_2$,
$\vert\be_3\vert=r_3$ and $r_1\in\{0,1,\ldots,r_2-1\}$, $r_2\in\{
1,\ldots,r_3-1\}$,
$r_3\in\{2,\ldots,N-1\}$.
\end{enumerate}
\end{enumerate}

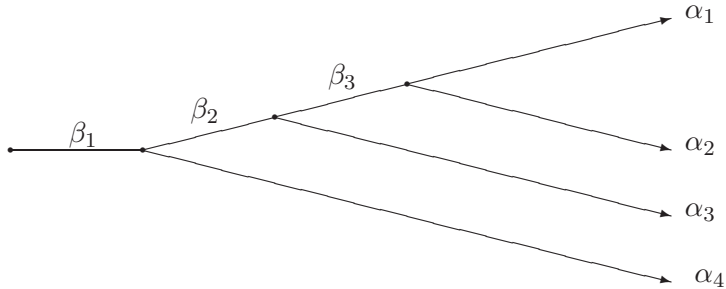
\begin{picture}(0,0)
\put(0,-50){\line(1,0){50}}
\put(50,-50){\line(4,1){50}}\put(50,-50){\vector(4,-1){200}}
\put(100,-37.5){\vector(4,1){150}}\put(100,-37.5){\vector(4,-1){150}}
\put(150,-25){\vector(4,-1){100}}\put(150,-25){\circle*{2}}
\put(0,-50){\circle*{2}}\put(50,-50){\circle*{2}}\put
(100,-37.5){\circle*{2}}
\put(255,0){$\al_1$}\put(255,-50){$\al_2$}\put(255,-75){$\al
_3$}\put(255,-100){
$\al_4$}\put(22,-47){$\be_1$}\put(68,-38){$\be_2$}\put
(120,-25){$\be_3$}
\end{picture}

\vspace{1.51in}Taking into consideration the possible resulting
diagrams, and defining $r(n)\in[0,r]$ by $r(n)=\frac{r}{n}$,
the following lemma can be shown.
\begin{lemma}
\label{lmmomentbp}
Given $\phi\in C^2_K(\R^{3\times d})$
and $\psi\in
C^2_K(\R^{4\times d})$, for all $n\in\N$, $t>0$, it follows that
%
%
\begin{eqnarray}\label{eqtmbpfinal}
\E\bigl\langle\phi,\bigl(\mu_t^{(n)}\bigr)^3\bigr\rangle&=&\bigl\langle Q^3_t\phi,\bigl(\mu
_0^{(n)}\bigr)^3\bigr\rangle+\frac{1}{n}\sum_
{r=0}^{
N-1}\sst{i}{j}^3
\bigl\langle\Ga^{(ij)}_{2;(r(n),t)}\phi,\bigl(\mu_0^ { (n) } \bigr)^2 \bigr\rangle
\nonumber\\[-12pt]\\[-12pt]
&&{} +
\frac{1}{n^2}\mathop{\sum_{r_1,r_2=0}}_{r_1<r_2}^{N-1}\sst
{i}{j}^3\bigl\langle\Ga^{ (12,ij)}_{1;(r_1(n),r_2(n),t)} ,\mu_0^{
(n)}\bigr\rangle+ o(1)\nonumber
\end{eqnarray}
and
%
%
\begin{eqnarray}\label{eqfmbpfinal}
\E\bigl\langle\psi,\bigl(\mu_t^{(n)}\bigr)^4\bigr\rangle&=&\bigl\langle Q^4_t\psi,\bigl(\mu
_0^{(n)}\bigr)^4\bigr\rangle+
\frac{1}{n}\sum_{r=0}^{N-1}\mathop{\sum_{i,j=1}}_{i\neq
j}^4\bigl\langle\Ga^{(ij)}_{3;(r(n),t)}\psi, \bigl(\mu_0^{(n)}\bigr)^3\bigr\rangle\nonumber\\[-4pt]
&&{}
+\frac{1}{n^2}\mathop{\sum_{r_1,r_2=0}}_{r_1<r_2}^{N-1}
\sst{i_1}{j_1}^4\sst{i_2}{j_2}
^3\bigl\langle\Ga^{(i_2j_2,i_1j_1)}_{2;(r_1(n),r_2(n),t)}\psi,\bigl(\mu
_0^{(n)}\bigr)^2\bigr\rangle\nonumber\\[-10pt]\\[-10pt]
&&{}+\frac{1}{n^3}\mathop{\mathop{\sum
_{r_k=0}}_{k=1,2,3}}_{r_1<r_2<r_3}^{N-1}\sst{i_1}{j_1}
^4\sst{i_2}{j_2}^3\bigl\langle\Ga
^{(12,i_2j_2,i_1j_1)}_{1;(r_1(n),r_2(n),r_3(n),t)} \psi, \mu_0^ { (n)
} \bigr\rangle
\nonumber\\[-4pt]
&&{} + o(1),\nonumber
\end{eqnarray}
where $\Ga^{\cdot}_{\cdot\cdot}$ is defined as in
(\ref{fnGamma}).\vspace*{-4pt}
\end{lemma}
\begin{pf}
See Appendix \ref{appA}.\vspace*{-4pt}
\end{pf}

Having now a formula for both the third and fourth moments of the branching
process,
with the exception of a some small technicalities to mention,
the moment formulae for the superprocess will follow almost immediately from
Lemmas \ref{lmEKconv3}, \ref{lmEKconv4} and \ref{lmmomentbp}.\vspace*{-4pt}
\begin{theorem}
\label{thmmomentsp}
Let $\phi\in C_K^\infty(\R^{3\times d})$ and $\psi\in
C^\infty_K(\R^{4\times d})$ be respectfully defined by
\[
\phi=\phi_1\otimes\phi_2\otimes\phi_3
\quad\mbox{and}\quad \psi=\psi_1\otimes\psi_2\otimes\psi_3\otimes\psi_4
\]
for $\phi_i$, $\psi_j$
$\in C^\infty_K(\R^{d})$, $i=1,2,3$, $j=1,2,3,4$. Then for all $0\leq
t_1\leq
t_2\leq t_3\leq t_4<\infty$,
%
%
\begin{eqnarray}\label{eqmixfmsp}
&&\E\langle\psi,\mu_{t_1}\mu_{t_2}\mu_{t_3}\mu_{t_4}\rangle\nonumber\\[-4pt]
&&\qquad=\bigl\langle\Ga^{ (0 , 0 , 0) } _ { 4;(t_1 , t_2,t_3,t_4)}\psi,\mu
_0^4\bigr\rangle+
\sst{i}{j}^4\int_0^{t_1}ds\bigl\langle\Ga
^{(ij,0,0,0)}_{3;(s,t_1,t_2,t_3,t_4) } \psi, \mu_0^3 \bigr\rangle
\nonumber\\
&&\qquad\quad{} +
\sst{i}{j}^3\int_{t_1}^{t_2} ds \bigl\langle\Ga
^{(0,ij,0,0)}_{3;(t_1,s,t_2,t_3,t_4)} \psi, \mu_0^3 \bigr\rangle+
\int_{t_2}^{t_3} ds \bigl\langle\Ga
^{(0,0,12,0)}_{3;(t_1,t_2,s,t_3,t_4)}\psi,\mu_0^3\bigr\rangle
\nonumber\\
&&\qquad\quad{} +
\sst{i_1}{j_1}^4\sst{i_2}{j_2}^3\int_0^{t_1} ds_2 \int_0^{s_2}
ds_1
\bigl\langle\Ga^{(i_2j_2,i_1j_1,0,0,0)}_{2;(s_1,s_2,t_1,t_2,t_3,t_4)}\psi
,\mu_0^2\bigr\rangle
\nonumber\\
&&\qquad\quad{}
+\sst{i_1}{j_1}^3\sst{i_2}{j_2}^3\int_{t_1}^{t_2} ds_2 \int_0^{t_1}
ds_1 \bigl\langle\Ga
^{(i_2j_2,0,i_1j_1,0,0)}_{2;(s_1,t_1,s_2,t_2,t_3,t_4)}\psi, \mu_0^2
\bigr\rangle\nonumber\\
&&\qquad\quad{} +
\sst{i}{j}^3\int_{t_1}^{t_2} ds_2 \int_{t_1}^{s_2} ds_1 \bigl\langle\Ga
^{(0,12,ij,0 ,0)}_{2;(t_1,s_1,s_2,t_2,t_3,t_4)}\psi, \mu_0^2 \bigr\rangle
\nonumber\\
&&\qquad\quad{} +
\sst{i}{j}^3\int_{t_2}^{t_3} ds_2 \int_0^{t_1} ds_1 \bigl\langle\Ga
^{(ij,0,0,12,0)} _{2;(s_1,t_1,t_2,s_2,t_3,t_4)}\psi, \mu_0^2 \bigr\rangle
\nonumber\\
&&\qquad\quad{} +
\int_{t_2}^{t_3} ds_2 \int_{t_1}^{t_2} ds_1 \bigl\langle\Ga
^{(0,12,0,12,0)}_{ 2;(t_1 ,s_1,t_2,s_2,t_3,t_4)}, \mu_0^2 \bigr\rangle
\nonumber\\
&&\qquad\quad{} +
\sst{i_1}{j_1}^4\sst{i_2}{j_2}^3\int_0^{t_1} ds_3 \int_0^{s_3}
ds_2\int_0^{
s_2} ds_1
\bigl\langle\Ga
^{(12,i_2j_2,i_1j_1,0,0,0)}_{1;(s_1,s_2,s_3,t_1,t_2,t_3,t_4)}\psi,\mu
_0\bigr\rangle
\\
&&\qquad\quad{} +
\sst{i_1}{j_1}^3\sst{i_2}{j_2}^3\int_{t_1}^{t_2} ds_3 \int_0^{t_1}
ds_2 \int_0^{s_2} ds_1 \bigl\langle\Ga
^{(12,i_2j_2,0,i_1j_1,0,0)}_{1;(s_1,s_2,t_1, s_3,t_2,t_3,t_4)}\psi,
\mu_0 \bigr\rangle\nonumber\\
&&\qquad\quad{} +
\sst{i}{j}^3\int_{t_1}^{t_2} ds_3 \int_{t_1}^{s_3} ds_2 \int_0^{t_1}
ds_1 \bigl\langle\Ga
^{(12,0,12,ij,0,0)}_{1;(s_1,t_1,s_2,s_3,t_2,t_3,t_4)}\psi,\mu
_0\bigr\rangle
\nonumber\\
&&\qquad\quad{} +
\sst{i}{j}^3\int_{t_2}^{t_3} ds_3 \int_0^{t_1} ds_2 \int_0^{s_2}
ds_1 \bigl\langle\Ga
^{(12,ij,0,0,12,0)}_{1;(s_1,s_2,t_1,t_2,s_3,t_3,t_4)}\psi,\mu
_0\bigr\rangle
\nonumber\\
&&\qquad\quad{} +
\int_{t_2}^{t_3} ds_3 \int_{t_1}^{t_2} ds_2 \int_0^{t_1} ds_1
\bigl\langle\Ga^{12 ,0,12,0,12,0)}_{1;(s_1,t_1,s_2,t_2,s_3,t_3,t_4)} \mu
_0 \bigr\rangle\nonumber
\end{eqnarray}
and
%
%
\begin{eqnarray}\label{eqmixtmsp}
\E\langle\phi,\mu_{t_1}\mu_{t_2}\mu_{t_3}\rangle
&=&\bigl\langle\Ga^{(0,0)}_{3;(t_1,t_2,t_3)}\phi,\mu_0^3\bigr\rangle+
\sst{i}{j}^3
\int_0^{t_1} ds \bigl\langle\Ga^{(ij,0,0)}_{2;(s,t_1,t_2,t_3)}\phi,\mu
_0^2\bigr\rangle\nonumber\\
&&{}
+\sst{i}{j}^3\int_0^{t_1} ds_2 \int_0^{s_2} ds_1 \bigl\langle\Ga
^{(12,ij,0,0)}_{ 1;(s_1,s_2,t_1,t_2,t_3)}\phi, \mu_0
\bigr\rangle\nonumber\\[-8pt]\\[-8pt]
&&{} +\int_{t_1}^{t_2} ds \bigl\langle\Ga^{(0,12,0)}_{2;(t_1,s,t_2,t_3)}\phi,
\mu_0^2\bigr\rangle\nonumber\\
&&{} +\int_{t_1}^{t_2} ds_2\int_0^{t_1} ds_1
\bigl\langle\Ga^{(12,0,12,0)}_{1;(s_1,t_1,s_2,t_2,t_3)}\phi,\mu
_0\bigr\rangle\nonumber.
\end{eqnarray}
\end{theorem}
\begin{pf}
See Appendix \ref{appB}.
\end{pf}

The purpose of the above moment formulae is due to the need for $L^2$ bounds,
the verification of which makes up the most essential part of this
paper. For
the remainder, any
arbitrary constant value, dependent only upon $0\leq T$, will be
denoted by
$C=C(T)$.
\begin{lemma}
\label{momentbound}Let $\phi\in S_d$, $d\leq3$, and define for
$x=(x_1,x_2,x_3)\in\R^{3\times d}$, $y=(y_1,y_2,y_3,y_4)\in\R
^{4\times d}$,
\[
\psi(x)\triangleq\phi(x_1-x_3)\phi(x_2-x_3)
\]
and
\[
\vp(y)\triangleq\phi(y_1-y_3)\phi(y_2-y_4).
\]
Suppose that $\mu=\{\mu_t\dvtx t\geq0\}$ is a superprocess over a
stochastic flow
such that $\mu_0\in M_F(\R^d)$ satisfies Assumption \ref{asn2}.
Then, for any $0\leq t_1\leq t_2\leq t_3\leq T<\infty$,
\[
\mbox{\textup{(i)}}\quad\int_0^T dt_3 \int_0^{t_3} dt_2 \int_0^{t_2} dt_1\, \E\langle\psi
,\mu_{t_1} \mu_{t_2}\mu_{t_3}\rangle\leq C \Vert\phi\Vert_{L^1}^2
\]
and
\[
\mbox{\textup{(ii)}}\quad\int_0^{t_3} {dt}_2 \int_0^{t_2} {dt}_1\, \E\langle\vp,\mu
_{t_1}\mu_{t_2}\mu_{ t_3}^2\rangle\leq C\Vert\phi\Vert_{L^1}^2.
\]
\end{lemma}
\begin{pf}
See Appendix \ref{appD}.
\end{pf}

We now proceed with the establishing existence of the GSILT.

\section{Existence of generalized self-intersection local time}
Generalized self-intersection local time (GSILT) at $u\in\R^d$, over
$B\subset\BA(\R^2)$, is defined formally as
\[
\LA(u;B)\triangleq\int_B dt\,ds \langle\delta_u,\mu_s\mu_t\rangle,
\]
where
$\delta_u(x)$
is the Dirac point-mass measure at $u$.

Note that in the above, and throughout the remainder of this paper, if
$\vp\dvtx\R^d\rightarrow\R$, the convention
\[
\langle\vp,\mu_s\mu_t\rangle=\int\mu_s(dx)\mu_t(dy)\vp(x-y)
\]
is made.

Since $\mu_s\mu_t=\mu_t\mu_s$, it makes sense to restrict GSILT
either above or
below the diagonal, and so we set
\[
\LA(u,T)=\LA\bigl(u;\{(s,t)\dvtx0\leq s\leq t\leq T\}\bigr)
\]
for fixed $T\in[0,\infty)$.

The above definition is clearly formal, and thus to make sense of this
a~limiting process will be constructed. For fixed $\lambda>0$, define
\[
G^{\lambda,u}(x)=\int_0^\infty dt\, e^{-\lambda t}q_t(u,x),
\]
then, since $\Gr$ is the resolvent to $L$ at $\lambda$, $(\lambda
-L)\Gr=\delta_u$ and
\mbox{$\Vert\Gr\Vert_{L^1}\leq\lambda^{-1}$}.

From Dynkin (\citeyear{dynkin65}), Theorem 0.5, it can be seen that
$\Gr(x)$ is not smooth (take, e.g., $x=u$), and thus it is
desired to estimate $\Gr$ by a class of smooth functions.

Since $\Gr\in L^1(\R^d)$, for any $\phi\in C_K^\infty(\R^d)$
\[
\langle\phi,\Gr\rangle\triangleq\int dx\, \Gr(x)\phi(x)<\infty,
\]
which implies $\Gr$ can be regarded as the element of $S'_d$ which sends
$\phi\in S_d$ to $\langle\phi,\Gr\rangle$. Thus, Lieb and Loss
(\citeyear{lieb}), Theorem 7.10, implies the
existence of a family $\{\G\dvtx\eps>0\}\subset C_K^\infty$ such that
$\G\rightarrow\Gr$ as
$\eps\rightarrow0$, in $S'_d$.

From H\"{o}rmander (\citeyear{horman}), $L$ is a continuous
operator on $S'_d$, and it is concluded
that
\[
\lim_{\eps\rightarrow0}(\lambda-L)\G=\delta_u,
\]
where convergence is in the sense of
distributions,
and so a limiting process is defined by
\[
\gamma^\lambda_\eps(u,T)\triangleq\int_0^T dt\int_0^t ds \langle
(\lambda-L)\G,\mu_s\mu_t\rangle,
\]
$\lambda>0$, $\eps>0$, $0\leq T<\infty$.

The goal now is to make sense of the operator $L$ appearing in the integrand.

\subsection{\texorpdfstring{An It\^{o} formula}{An Ito formula}}

As in the independent case, the derivation of the evolution equation is
accomplished through the construction, and careful application, of an
appropriate It\^{o} formula. This construction will mimic that of
Adler and Lewin (\citeyear{adler}), which begins with application of
It\^{o}'s lemma
to the nonanticipative
functional $f$, given by
\[
f(t,x)=x\int_0^tds\, \mu_s(\psi),
\]
where $\psi\in C^2_K(\R^d)$,\vspace*{1pt} and $x$ is a
$\R$-valued random variable. Note that from the SPDE (\ref{eqSPDE}), if
$\phi\in
C^\infty_K(\R^d)$, then $\mu_t(\phi)$ is a continuous
semi-martingale with
decomposition
\[
\mu_t(\phi)=\mu_0(\phi)+Z_t(\phi)+V_t(\phi),
\]
where
\[
V_t(\phi)\triangleq\int_0^t d{s}\,\mu_s(L\phi).
\]
\begin{theorem}
\label{thmsemimartingale}
If $\phi\in S_d$ then $\mu_t(\phi)$ is an a.s.
continuous semimartingale.
\end{theorem}
\begin{pf}
See Appendix \ref{appD}.
\end{pf}

Through some careful work (outlined in the \hyperref[app]{Appendix}),
we arrive at the
following.
\begin{lemma}
\label{ito}Given $\Psi\in S_{2d}$,
\begin{eqnarray*}
\int_0^T dt \int_0^t ds \langle L_2\Psi,\mu_s\mu_t\rangle
& = &
\int_0^Tdt \langle\Psi,\mu_t\mu_T\rangle-\int_0^T dt \langle
\Psi,\mu_t\mu_t\rangle
\\
&&{} -\int_0^T\int_{\R^d} Z(dt,dy)\int_0^t ds \langle\Psi(\cdot
,y),\mu_s\rangle,
\end{eqnarray*}
where
\[
(L_2\Psi)(x,y)\triangleq\frac{1}{2}\sum_{i,j=1}^d
a_{ij}(y)\,\pa_{2_i}\,\pa_{2_j}\Psi(x,y),
\]
and $Z(dt,dy)$ is the corresponding martingale measure.
\end{lemma}
\begin{pf}
See Appendix \ref{appE}.
\end{pf}\vfill\eject

\subsection{Existence}

Using lemma \ref{ito} with $\G$ in place of $\Psi$, we now have
\begin{eqnarray*}
\gamma^\lambda_\eps(u,T)&=&\lambda\int_0^T dt\int_0^s ds \langle
\G,\mu_s\mu_t\rangle\\[-2pt]
&&{}-\int_0^T dt \langle\G,\mu_t\mu_T\rangle
+ \int_0^T dt \langle\G,\mu_t\mu_t\rangle\\[-2pt]
&&{}+ \int_0^T \int_{\R
^d}
Z(dt,dy) \int_0^t ds \langle\G(\cdot-y),\mu_s\rangle.
\end{eqnarray*}

As in Rosen (\citeyear{rosen}) and Adler and Lewin (\citeyear
{adler,lewin}), the
issue of ``local double points'' must be addressed, that is, the set of points
lying on the diagonal in~$\R^2$, which will be (falsely) counted as
points of
self-intersection when $u=0$, and will lead to singularities in dimensions
greater than one. Due to this we follow the idea first proposed by
Adler and
Lewin, and renormalize our GSILT via subtraction of the term involving ``local
double points.'' It is easy enough to see that the term involving the ``local
double points'' is given by
$\int_0^T dt \langle\G,\mu_t\mu_t\rangle$, and thus we define our
renormalized
limiting process to generalized self-intersection local time at $u\in
\R^d$, over
the set $\{(s,t)\dvtx0\leq s< t\leq T\}$ by
\begin{eqnarray*}
\LA_\eps^\lambda(u,T)
& = &\gamma_\eps^\lambda(u,T)-\int_0^T dt \langle\G,\mu_t\mu
_t\rangle\\[-2pt]
& = &\lambda\int_0^T dt\int_0^t ds \langle\G,\mu_s\mu_t\rangle-
\int_0^T dt \langle\G,\mu_t\mu_T\rangle\\[-2pt]
&&{} + \int_0^T \int_{\R^d} Z(dt,dy) \int_0^t ds \langle\G(\cdot
-y),\mu_s\rangle.
\end{eqnarray*}

Using Lemma \ref{momentbound} existence follows almost immediately.
\begin{theorem}
\label{thmexistence}
Suppose that $\mu=\{\mu_t\dvtx t\geq0\}$ is a
$d$-dimensional superprocess over a stochastic flow such that $\mu
_0\in
M_F(\R^d)$, $d\leq3$, satisfies Assumption~\ref{asn2}.
Fix $T\in[0,\infty)$ and define $\L$ as above, then for $0\leq
s<t\leq T$,
\[
L^2 - \lim_{\eps\rightarrow0}\L=\LA^\lambda(u,T),
\]
uniformly in $u\in\R^d$, where $\LA^\lambda(u,T)$ is defined by
\begin{eqnarray*}
\LA^\lambda(u,T)& = &\lambda\int_0^T dt\int_0^t ds \langle\Gr,\mu
_s\mu_t\rangle-
\int_0^T dt \langle\Gr,\mu_t\mu_T\rangle\\[-2pt]
&&{} + \int_0^T \int_{\R^d}
Z(dt,dy) \int_0^t ds \langle\Gr(\cdot-y),\mu_s\rangle.
\end{eqnarray*}
For each $\lambda>0$, $\LA^\lambda(u,T)$ is referred to as the
self-intersection
local time at~$u$, up to time $T$, for a superprocess over a stochastic
flow.\vadjust{\goodbreak}
\end{theorem}
\begin{pf}
See Appendix \ref{appF}.
\end{pf}

It should be noted that, as with the the SILT of Adler and Lewin
(\citeyear{adler}), for $d>3$ the GSILT can be shown to blow up to
infinity. It remains an open question if renormalization processes,
such as those of Rosen (\citeyear{rosen}), exist for dimensions $d>3$.

\begin{appendix}\label{app}
\section{\texorpdfstring{Proof of Lemma \lowercase{\protect\ref{lmmomentbp}}}{Proof of Lemma 3.8}}\label{appA}
We now proceed with the moment calculations. Much of what
follows will be a consequence of the Markov property, and the reader is referred
to Skoulakis and Adler (\citeyear{skoulakis}) for a similar
calculation for the first
and second moments. Note that if $t\geq0$ and $r\in\N$, we define
$N\in\N$ and
$r(n)\in[0,r]$ by $N=[nt]$ and $r(n)=\frac{r}{n}$.
Recall,
%
%
\begin{equation}
\E\bigl\langle\phi,\bigl(\mu_t^{(n)}\bigr)^4\bigr\rangle=\frac{1}{n^4} \mathop{\sum
_{\al_j\thicksim_n
t}}_{j=1,2,3,4} \E\phi(Y^{\al_1,n}_t,Y^{\al_2,n}_t,Y^{\al_3,n}_t,Y^{
\al_4,n}_t)\E\prod_{ i=1 } ^4 1_{ \al_i,n}(t).\hspace*{-24pt}
\end{equation}
If $\al_1(0),\al_2(0),\al_3(0)$, and $\al_4(0)$
are given, case (I) gives
\[
\E\phi(Y^{\al_1,n}_t,Y^{\al_2,n}_t,Y^{\al_3,n}_t,Y^{\al_4,n}_t) =
Q^4_t\phi\bigl(x_{\al_1(0)},x_{\al_2(0)},x_{\al_3(0)},x_{\al_4(0)}\bigr)
\]
and $\E\prod_{i=1}^41_{\al_i,n}(t) = (\frac{1}{2})^{4N}$.
For any $\al_1(0)$, $\al_2(0)$, $\al_3(0)$, $\al_4(0)$, there are~$2^{4N}$
corresponding $(\al_1,\al_2,\al_3,\al_4)$ which result from binary branching
over~$N$ steps. We thus arrive at the following contribution from case (I):
\begin{eqnarray*}
&&\frac{1}{n^4} \mathop{\sum_{\al_k(0)=1, k=1,2,3,4}}_{\al
_\ell(0)\neq\al_k(0), \ell\neq k}
Q^4_t\phi\bigl(x_{\al_1(0)},x_{\al_2(0)},x_{\al_3(0)},x_{\al
_4(0)}\bigr)\nonumber\\
&&\qquad =\frac{1}{n^4}\sum_{\al_k(0)=1,k=1,2,3,4}
Q^4_t\phi\bigl(x_{\al_1(0)},x_{\al_2(0)},x_{\al_3(0)},x_{\al
_4(0)}\bigr)\nonumber\\
&&\qquad\quad{} -\frac{1}{n^4}\mathop{\sum_{\al_k(0)=1}}_{k=1,2,3}\sst{i}{j}^4
(\Phi_{ij}Q^4_t\phi)\bigl(x_{\al_1(0)},x_{\al_2(0)},x_{\al
_3(0)}\bigr)\nonumber\\
&&\qquad\quad{} -\frac{1}{n^4}\mathop{\sum_{\al_k(0)=1}}_{j=1,2}\sst
{i_1}{j_1}^4\sst{i_2}{j_2}^3
(\Phi_{i_2j_2}\Phi_{i_1j_1}Q^4_t\phi)\bigl(x_{\al_1(0)},x_{\al_2(0)}\bigr)\\
&&\qquad\quad{}-\frac{1}{n^4}\sum_{
\al_1(0)=1}\sst{i_1}{j_1}^4\sst{i_2}{j_2}^3
(\Phi_{12}\Phi_{i_2j_2}\Phi_{i_1j_1}Q^4_t\phi)\bigl(x_{\al_1(0)}\bigr),
\end{eqnarray*}
which by the definition of $\mu^{(n)}$,
\begin{eqnarray*}
&&\frac{1}{n^4} \mathop{\sum_{\al_k(0)=1, k=1,2,3,4}}_{
\al_\ell(0)\neq\al_k(0), \ell\neq k}
Q^4_t\phi\bigl(x_{\al_1(0)},x_{\al_2(0)},x_{\al_3(0)},x_{\al_4(0)}\bigr)
\\[-2pt]
&&\qquad
=\bigl\langle Q^4\phi,\bigl(\mu_0^{(n)}\bigr)^4\bigr\rangle-
\frac{1}{n}\bigl\langle\Phi_{ij}Q^4_t\phi,\bigl(\mu_0^{(n)}\bigr)^3\bigr\rangle\\[-2pt]
&&\qquad\quad{} - \frac{1}{n^2}\bigl\langle\Phi_{i_2j_2}\Phi_{i_1j_1}Q^4_t\phi,\bigl(\mu
_0^{(n)}\bigr)^2\bigr\rangle-
\frac{1}{n^3}\bigl\langle\Phi_{12}\Phi_{i_2j_2}\Phi_{i_1j_1}Q^4_t\phi
,\mu_0^{(n)}\bigr\rangle.
\end{eqnarray*}
From Lemma \ref{lmEKconv2} all but the first term on the right-hand side
will vanish as $n\rightarrow\infty$, and thus
%
%
\begin{eqnarray}\label{eqfmbpI}
&&\frac{1}{n^4} \mathop{\sum_{\al_k(0)=1, k=1,2,3,4}}_{
\al_\ell(0)\neq\al_k(0), \ell\neq k}
Q^4_t\phi\bigl(x_{\al_1(0)},x_{\al_2(0)},x_{\al_3(0)},x_{\al_4(0)}\bigr) \nonumber\\[-9pt]\\[-9pt]
&&\qquad=
\bigl\langle Q^4_t\phi,\bigl(\mu_0^{(n)}\bigr)^4\bigr\rangle+
o(1).\nonumber
\end{eqnarray}

For case (II), given $\al_1(0),\al_2(0),\be(0)$ and $r$, proceeding
as before,
\begin{eqnarray*}
&&\E\phi(Y^{\al_1,n}_t,Y^{\al_2,n}_t,Y^{\al_3,n}_t,Y^{\al_4,n}_t)\\[-2pt]
&&\qquad=\frac{1}{12}\sst{i}{j}^4
\bigl(Q^3_{r(n)}\Phi_{ij}Q^4_{t-r(n)}\phi\bigr)\bigl(x_{\al_1(0)},x_{\al
_2(0)},x_{\be(0)}\bigr),
\end{eqnarray*}
and if for any distinct $i,j\in\{1,2,3,4\}$ we define $i',j'$ to be the
exhaustive elements of $\{1,2,3,4\}\setminus\{i,j\}$,
\begin{eqnarray*}
\E\prod_{i=1}^41_{\al_i,n}(t)&=&(\E1_{\al_{i'},n}(t)
)(\E1_{\al_{
j'},n}(t))\E\E\bigl[1_{\al_i,n}(t)1_{\al_j,n}(t){|}
\FA^n_{r(n)}\bigr]\\[-2pt]
& = &2^{-(4N-r-1)}.
\end{eqnarray*}
For any $\al_1(0)$, $\al_2(0)$, $\be_1(0)$ and $r$, there are $2^{4N-r-1}$
corresponding\vspace*{1pt} tuples $(\al_1,\al_2,\al_3,\al_4)$ which result from binary
branching over $N$ steps and $2\cdot\binom{4}{2}$ possible
arrangements for
$(\al_1,\al_2,\al_3,\al_4)$. We thus arrive at the following
contribution from
case (II):
\begin{eqnarray*}
&&\frac{1}{n^4}
\sum_{r=0}^{N-1}\sst{i}{j}^4\mathop{\mathop{\sum_{\al_1(0),\al
_2(0),\be(0)=1}}_{\al_1(0)\neq\al_2(0),\al_\ell(0)\neq\be
(0)}}_{\ell=1,2}
\bigl(\Ga^{(ij)}_{3;(r(n),t)}\phi\bigr)\bigl(x_{\al_1(0)},x_{\al_2(0)},x_{\be
(0)}\bigr)\\[-2pt]
&&\qquad =\frac{1}{n}\sum_{r=0}^{N-1}\sst{i}{j}^4\bigl\langle\Ga
^{(ij)}_{3;(r(n),t)}\phi, \bigl(\mu_0^{(n)}\bigr)^3\bigr\rangle\\[-2pt]
&&\qquad\quad{} -\frac{1}{n^2}\sum_{r=0}^{N-1}\sst{i_1}{j_1}^4\sst{i_2}{j_2}^3
\bigl\langle\Phi_{i_2j_2}\Ga^{(i_1j_1)}_{3;(r(n),t)}\phi,\bigl(\mu
_0^{(n)}\bigr)^2\bigr\rangle\\[-2pt]
&&\qquad\quad{} -\frac{1}{n^3}\sum_{r=0}^{N-1}\sst{i_1}{j_1}^4\sst{i_2}{j_2}^3
\bigl\langle\Phi_{12}\Phi_{i_2j_2}\Ga^{(i_1j_1)}_{3;(r(n),t)}\phi,\bigl(\mu
_0^{(n)}\bigr) \bigr\rangle.
\end{eqnarray*}
Again from Lemma \ref{lmEKconv3}, all but the first term on the
right-hand side
will vanish as $n\rightarrow\infty$ and thus,
%
%
\begin{eqnarray}\label{eqfmbpII}
&&\frac{1}{n^4}
\sum_{r=0}^{N-1}\sst{i}{j}^4\mathop{\mathop{\sum_{\al_1(0),\al
_2(0),\be(0)=1}}_{
\al_1(0)\neq\al_2(0),\al_\ell(0)\neq\be(0)}}_{\ell=1,2}
\bigl(\Ga^{(ij)}_{3;(r(n),t)}\phi\bigr)\bigl(x_{\al_1(0)},x_{\al_2(0)},x_{\be
(0)}\bigr)\nonumber\\[-8pt]\\[-8pt]
&&\qquad=\frac{1}{n}\sum_{r=0}^{N-1}\sst{i}{j}^4\bigl\langle\Ga
^{(ij)}_{3;(r(n),t)}\phi, \bigl(\mu_0^{(n)}\bigr)^3\bigr\rangle+ o(1).\nonumber
\end{eqnarray}

Cases (III) and (IV) will now be considered together.
For case (III), given $\be_1(0)$, $\be_2(0)$, $r_1$ and $r_2$,
\begin{eqnarray*}
&&\E\phi(Y^{\al_1,n}_t,Y^{\al_2,n}_t,Y^{\al_3,n}_t,Y^{\al_4,n}_t)\\
&&\qquad =\frac{1}{12}\sst{i_1}{j_1}^4\mathop{\mathop{\sum
_{i_2,j_2=1}}_{i_2\neq
j_2}}_{i_2,j_2\neq
i_1}^3\bigl(\Ga^{(i_2j_2,i_1j_1)}_{2;(r_1(n),r_2(n),t)}
\phi\bigr)\bigl(x_{\be_1(0)},x_{\be_2(0)}\bigr)
\end{eqnarray*}
and
\[
\E\prod_{i=1}^41_{\al_i,n}(t) = 2^{-(4N-r_1-r_2-2)}.
\]

For case (IV), given $\al(0)$, $\be_1(0)$, $r_1$ and $r_2$,
\begin{eqnarray*}
&&\E\phi(Y^{\al_1,n}_t,Y^{\al_2,n}_t,Y^{\al_3,n}_t,Y^{\al_4,n}_t)\\
&&\qquad =\frac{1}{48}\sst{i_1}{j_1}^4\mathop{\mathop{\sum
_{i_2,j_2=1}}_{i_2\neq j_2}}_{i_2= i_1
\ \mathrm{or}\ j_2=
i_1}^3\bigl(\Ga^{(i_2j_2,i_1j_1)}_{2;(r_1(n),r_2(n),t)}\phi\bigr)
\bigl(x_{\al
(0)},x_{\be_1(0)}
\bigr)
\end{eqnarray*}
and
\[
\E\prod_{i=1}^41_{\al_i,n}(t) = 2^{-(4N-r_1-r_2-2)}.
\]

Given two initial ancestors, there are $2^{4N-r_1-r_2-2}$ possible
trees, and a~possible $2\cdot\binom{4}{2}$ arrangements for $\al_1$, $\al_2$,
$\al_3$,
$\al_4$ upon each tree (requiring $r_1<r_2$) that result in case (III).
Furthermore, there are $2^{4N-r_1-r_2-2}$ possible trees, and a possible
$2\cdot\binom{2}{1}\cdot\binom{3}{2}\cdot\binom{4}{3}$
arrangements for
$\al_1$, $\al_2$, $\al_3$, $\al_4$ upon each tree that result in
case (IV).\vspace*{1pt} It
follows that the contribution coming from the sum of case (III) and
case (IV) is
given by
\begin{eqnarray*}
&&\frac{1}{n^4}\mathop{\sum_{r_1,r_2=0}}_{r_1<r_2}^{N-1}\sst{i_1}{j_1}
^4\sst{ i_2 } { j_2 }
^3 \mathop{\mathop{\sum_{\al(0)=1}}_{\be(0)=1}}_{\al(0)\neq
\be(0)}\bigl(\Ga^{(i_2j_2,
i_1j_1)}_{2;(r_1(n),r_2(n),t)}\phi\bigr)\bigl(x_{\al(0)},x_{\be(0)}
\bigr)\\
&&\qquad
=\frac{1}{n^2}\mathop{\sum_{r_1,r_2=0}}_{r_1<r_2}^{N-1}\sst
{i_1}{j_1}^4\sst{
i_2 } { j_2 }
^3\bigl\langle\Ga^{(i_2j_2,i_1j_1)}_{2;(r_1(n),r_2(n),t)}\phi,\bigl(\mu
_0^{(n)}\bigr)^2\bigr\rangle\\
&&\qquad\quad{} -\frac{1}{n^3}\mathop{\sum_{r_1,r_2=0}}_{r_1<r_2}^{N-1}\sst
{i_1}{j_1}^4\sst{
i_2 } { j_2 }
^3\bigl\langle\Phi_{12}\Ga^{(i_2j_2,i_1j_1)}_{2;(r_1(n),r_2(n),t)}\phi
,\mu_0^{(n)}\bigr\rangle.
\end{eqnarray*}
Thus, again from Lemma \ref{lmEKconv3}, the second term vanishes as
$n\rightarrow\infty$, and we have the contribution
%
%
\begin{equation}\label{eqfmbpIII}\qquad
\frac{1}{n^2}\mathop{\mathop{\sum
_{r_1=0}}_{r_2=0}}_{r_1<r_2}^{N-1}\sst{i_1}{j_1}^4\sst{
i_2}{j_2} ^3\bigl\langle\bigl\langle\Phi_{12}\Ga
^{(i_2j_2,i_1j_1)}_{2;(r_1(n),r_2(n),t)}\phi,\mu_0^{(n)}\bigr\rangle
,\bigl(\mu_0^{ (n)}\bigr)^2\bigr\rangle+ o(1).\
\end{equation}

Considering subcase (V)(A), given $r_1$, $r_2$, $r_3$ and $\be_1(0)$,
\begin{eqnarray*}
&&\E\phi(Y^{\al_1,n}_t,Y^{\al_2,n}_t,Y^{\al_3,n}_t,Y^{\al_4,n}_t)\\
&&\qquad=
\frac{1}{6}\sst{i}{j}^4\E\E\bigl[\bigl(\Phi_{ij}Q^4_{t-r_3(n)}\phi
\bigr)\bigl(Y^{\al_1,n}_{
r_3(n)},Y^{\al_2,n}_{r_3(n)},Y^{\be_3,n}_{r_3(n)}\bigr){|}\FA
^n_{r_2(n)}\bigr]
\\
&&\qquad= \frac{1}{12} \sst{i_1}{j_1}^4\mathop{\mathop{\sum
_{i_2,j_2=1}}_{i_2\neq
j_2}}_{i_2,j_2\neq i_1}^3
\bigl(\Ga^{(12,i_2j_2,i_1j_1)}_{1;(r_1(n),r_2(n),r_3(n),t)}\phi\bigr)\bigl(x_{\be_1(0)}\bigr).
\end{eqnarray*}
Furthermore,
\[
\E\prod_{i=1}^41_{\al_i,n}(t)=2^{-(4N-r_3-r_2-r_1-3)}.
\]

For subcase (V)(B), given $r_1$, $r_2$, $r_3$ and $\be_1(0)$,
\begin{eqnarray*}
&&\E\phi(Y^{\al_1,n}_t,Y^{\al_2,n}_t,Y^{\al_3,n}_t,Y^{\al_4,n}_t)\\
&&\qquad =\frac{1}{12}\sst{i}{j}^4\E\E\bigl[\bigl(\Phi_{ij}Q^4_{t-r_3(n)}\phi
\bigr)\bigl(Y^{\al_1,n}_{
r_3(n)},Y^{\al_2,n}_{r_3(n)},Y^{\be_3,n}_{r_3(n)}\bigr){|}\FA
^n_{r_2(n)}\bigr]
\\
&&\qquad=\frac{1}{48} \sst{i_1}{j_1}^4 \mathop{\mathop{\sum
_{i_2,j_2=1}}_{i_2\neq
j_2}}_{i_2=i_1\ \mathrm{or}\ j_2= i_1}^3
\bigl(\Ga^{(12,i_2j_2,i_1j_1)}_{1;(r_1(n),r_2(n),r_3(n),t)}\phi\bigr)\bigl(x_{\be_1(0)}\bigr)
\end{eqnarray*}
and
\[
\E\prod_{i=1}^41_{\al_i,n}(t) = 2^{-(4N-r_1-r_2-r_3-3)}.
\]

Given one initial ancestor there are $2^{4N-r_1-r_2-r_3-3}$ possible
trees, and
a~possible $\binom{1}{1}\cdot\binom{2}{1}\cdot\binom{4}{2}$
arrangements for
$\al_1$, $\al_2$, $\al_3$, $\al_4$ upon each tree (requiring
$r_2<r_3$) that
result in case (V)(A). Furthermore, there are $2^{4N-r_1-r_2-r_3-3}$ possible
trees, and a possible
$\binom{1}{1}\cdot\binom{2}{1}\cdot\binom{3}{2}\cdot\binom{4}{3}$
arrangements\vspace*{1pt} for $\al_1$, $\al_2$, $\al_3$, $\al_4$ upon each tree
that result
in case (V)(B). It follows that the contribution coming from the sum of
subcase (V)(A) and subcase (V)(B) is given by
%
%
\begin{equation}\label{eqfmbpIV}
\frac{1}{n^3}\mathop{\sum_{r_1,r_2,r_3=0}}_{r_1<r_2<r_3}^{N-1}\sst{i_1}{j_1}
^4\sst{i_2}{j_2}^3\bigl\langle\Ga^{(12,i_2j_2,i_1j_1)}_{1;(r_1(n),r_2(n),
r_3(n),t)}\phi, \mu_0^{(n)}\bigr\rangle.
\end{equation}

Therefore, from (\ref{eqfmbpI}), (\ref{eqfmbpII}),
(\ref{eqfmbpIII}) and (\ref{eqfmbpIV}), the lemma is shown.

\section{\texorpdfstring{Proof of Theorem \lowercase{\protect\ref{thmmomentsp}}}{Proof of Theorem 3.9}}\label{appB}

Again, due to similarity and escalating difficulty, we forgo the proof
of the
third moment in favor of the fourth moment. We first prove a needed lemma.
\begin{lemma}
Given $\phi_k, \psi_j\in C^\infty_K(\R^{d})$, $k=1,2,3$, $j=1,2,3,4$,
let $\phi=\phi_1\otimes\phi_2\otimes\phi_3$ and
$\psi=\psi_1\otimes\psi_2\otimes\psi_3\otimes\psi_4$. For
any $t\geq0$, the following hold:
%
%
\begin{eqnarray}\label{eqtmsp}
\E\langle\phi,\mu_t^3\rangle
&=& \langle Q^3_t\phi,\mu_0^3\rangle+
\sst{i}{j}^3 \int_0^t ds
\bigl\langle\Ga^{(ij)}_{2;(s,t)}\phi,\mu_0^2\bigr\rangle\nonumber\\[-8pt]\\[-8pt]
&&{} + \sst{i}{j}^3\int
_0^t ds_2
\int_0^{s_2} ds_1
\bigl\langle\Ga^{(12,ij)}_{1;(s_1,s_2,t)}\phi,\mu_0\bigr\rangle\nonumber
\end{eqnarray}
and
%
%
\begin{eqnarray}\label{eqfmsp}\quad
\E\langle\psi,\mu_t^4\rangle&=&\langle Q^4_t\psi,\mu_0^4\rangle+
\sst{i}{j}^4\int_0^t ds \bigl\langle\Ga^{(ij)}_{3;(s,t)}\psi,\mu
_0^3\bigr\rangle
\nonumber\\[-3pt]
&&{} +\sst{i_1}{j_1}^4\sst{i_2}{j_2}^3
\int_0^t ds_2\int_0^{s_2} ds_1\bigl\langle\Ga
^{(i_2j_2,i_1j_1)}_{2;(s_1,s_2,t)}\psi, \mu_0^2 \bigr\rangle\\[-3pt]
&&{} + \sst{i_1}{j_1}^4\sst{i_2}{j_2}^3\int_0^t ds_3\int_0^{s_3} ds_2
\int_0^{s_2} ds_1\bigl\langle\Ga
^{(12,i_2j_2,i_1j_1)}_{1;(s_1,s_2,s_3,t)}\psi,\mu_0\bigr\rangle.
\nonumber\vspace*{-3pt}
\end{eqnarray}
\end{lemma}
\begin{pf}
To begin, note that Lemma \ref{lmEKconv2} implies
$(\mu_0^{(n)})^\ell\Rt\mu_0^\ell$ for any $\ell\in\N$, and thus
the first term of the
right-hand sides of (\ref{eqtmbpfinal}) and (\ref{eqfmbpfinal})
converge,
respectively, to the first term of the right-hand sides of (\ref
{eqtmsp}) and
(\ref{eqfmsp}) as $n\rightarrow\infty$. Since $Q^k_t$ is a strongly
continuous
contraction semigroup for \mbox{$k\in\N$} (Lemma \ref{lmtf}), for any $\phi
\in
C_b(\R^d)$ which satisfies $\Vert\phi\Vert_{bL}=1$, \mbox{$\Vert
Q^k_t\phi\Vert_\infty\leq1$},
and $\sup_{x\neq y}\frac{\vert Q^k_t\phi(x)-Q^k_t\phi(y)\vert
}{\vert x-y\vert}\leq
1$. Thus, for any $k\in N$, $\Vert\phi\Vert_{bL}=1$ implies
$\Vert Q^k_t\phi\Vert_{bL}\leq1$. From Lemma \ref{lmEKconv3}, the
remaining terms
on the right-hand sides of (\ref{eqtmbpfinal}) and~(\ref{eqfmbpfinal})
converge,
respectively, to the remaining terms of the right-hand sides of~(\ref{eqtmsp})
and~(\ref{eqfmsp}) as $n\rightarrow\infty$. It remains to show that the
left-hand sides of~(\ref{eqtmbpfinal}) and~(\ref{eqfmbpfinal}) converge, respectively, to
the left-hand
sides of (\ref{eqtmsp}) and~(\ref{eqfmsp}), but this follows
immediately from
Lemma \ref{lmEKconv4}\vspace*{-3pt}
\end{pf}

The proof of the main theorem can now be shown.\vspace*{-3pt}
\begin{pf*}{Proof of Theorem \ref{thmmomentsp}}
Using the Markov property and Lemma \ref{lmmoment12sp}, it follows that
\begin{eqnarray*}
&&\E\langle\psi,\mu_{t_1}\mu_{t_2}\mu_{ t_3}\mu_{t_4}\rangle\\[-2pt]
&&\qquad=
\E\mu_{t_1}(\psi_1)\mu_{t_2}(\psi_2)\mu_{t_3}^2(\psi_3\otimes
Q_{t_4-t_3}\psi_4)\\[-2pt]
&&\qquad= \E\mu_{t_1}(\psi_1)\mu_{t_2}^3\bigl(\psi_2\otimes Q^2_{t_3-t_2}(\psi
_3\otimes
Q_{t_4-t_3}\psi_4)\bigr)\\[-2pt]
&&\qquad\quad{} +
\int_{0}^{t_3-t_2} ds \,\mu_{t_1}(\psi_1)\\[-2pt]
&&\qquad\quad\hspace*{40.2pt}{}\times\mu_{t_2}^2\bigl(\psi_2\otimes
Q_{s}\Phi_{12}Q^2_{t_3-t_2-s}(\psi_3\otimes
Q_{t_4-t_3}\psi_4)\bigr)\\[-2pt]
&&\qquad= \E\mu_{t_1}(\psi_1)\mu_{t_2}^3\bigl(\pi_1 Q^2_{t_3-t_2}\pi_1
Q_{t_4-t_3}(\psi_2\otimes\psi_3\otimes\psi_4)\bigr)\\[-2pt]
&&\qquad\quad{} + \int_{t_2}^{t_3} ds\, \E\mu_{t_1}(\psi_1)\\[-2pt]
&&\qquad\quad\hspace*{27.8pt}{}\times\mu_{t_2}^2\bigl(\pi_1
Q_{s-t_2}\Phi_{12}Q^2_{t_3-s}\pi_1Q_{t_4-t_3}(\psi_2\otimes\psi
_3\otimes
\psi_4)\bigr)\\[-2pt]
&&\qquad=\E\mu_{t_1}^4\bigl(\pi_1\Ga^{(0,0)}_{3;(t_2-t_1,t_3-t_1,t_4-t_1)}\psi
\bigr)\\[-2pt]
&&\qquad\quad{} + \sst{i}{j} ^3\int_{t_1}^{t_2} ds\, \E\mu_{t_1}^3\bigl(\pi_1\Ga^{(ij,0,0)}_{
2;(s-t_1,t_2-t_1,t_3-t_1,t_4-t_1)}\psi\bigr)\\[-2pt]
&&\qquad\quad{}
+\int_{t_2}^{t_3} ds\, \E\mu_{t_1}^3\bigl(\pi_1\Ga^{(0,12,0)}_{2;(t_2-t_1,
s-t_1,t_3-t_1,t_4-t_1)}\psi\bigr)\\[-2pt]
&&\qquad\quad{}
+\sst{i}{j}^3\int_{t_1}^{t_2} ds_2 \int_{t_1}^{s_2} ds_1\, \E
\mu_{t_1}^2\bigl(\pi_1\Ga^{(12,ij,0,0)}_{1;(s_1-t_1,s_2-t_1,t_2-t_1,t_3-t_1,
t_4-t_1)}\psi\bigr)\\[-2pt]
&&\qquad\quad{}
+\int_{t_2}^{t_3} ds_2 \int_{t_1}^{t_2} ds_1 \,\E\mu_{t_1}
^2\bigl(\pi_1\Ga
^{(12,0,12,0)}_{1;(s_1-t_1,t_2-t_1,s_2-t_1,t_3-t_1,t_4-t_1)}\psi\bigr).
\end{eqnarray*}

To make sense of the remainder of the proof, each of the above five
terms will
now be considered separately.

From (\ref{eqfmsp}),
%
%
\begin{eqnarray}\label{eqfmp1}
\qquad&&\E\mu_{t_1}^4\bigl(\pi_1\Ga^{(0,0)}_{3;(t_2-t_1,t_3-t_1,t_4-t_1)}\psi
\bigr)\nonumber\\
&&\qquad = \mu_{0}^4 \bigl(\Ga^{(0,0,0)}_{4;(t_1,t_2,t_3,t_4)}\psi\bigr)
+
\sst{i}{j}^4\int_0^{t_1} ds \,\mu_{0}^3\bigl(
\Ga^{(ij,0,0,0)}_{3;(s,t_1,t_2,t_3,t_4)}\psi\bigr)\\
&&\qquad\quad{} +\sst{i_1}{j_1}^4\sst{i_2}{j_2}^3
\int_0^{t_1} ds_2\int_0^{s_2} ds_1 \,\mu_{0}^2\bigl(\Ga
^{(i_2j_2,i_1j_1,0,0,0)}_{
2;(s_1,s_2,t_1 , t_2 , t_3 , t_4)}\psi\bigr)\nonumber\\
&&\qquad\quad{} + \sst{i_1}{j_1}^4\sst{i_2}{j_2}^3\int_0^{t_1} ds_3\int_0^{s_3} ds_2
\int_0^{s_2} ds_1 \,\mu_{0}\bigl(\Ga
^{(12,i_2j_2,i_1j_1,0,0,0)}_{1;(s_1,s_2,s_3,t_1,
t_2 , t_3 , t_4) } \psi\bigr).\nonumber
\end{eqnarray}

From (\ref{eqtmsp}),
%
%
\begin{eqnarray}\label{eqfmp2}
\quad&&\sst{i}{j}^3 \int_{t_1}^{t_2} ds\,\E\mu_{t_1}^3 \bigl(\pi_1\Ga
^{(ij,0,0)
} _ {
2;(s-t_1,t_2-t_1,t_3-t_1,t_4-t_1)}\psi\bigr)\nonumber\\
&&\qquad = \sst{i}{j}^3 \int
_{t_1}^{t_2}
ds \,\mu_0^3\bigl(\Ga^{(0,ij,0,0)}_{
3;(t_1,s,t_2,t_3,t_4)}\psi\bigr)\nonumber\\[-8pt]\\[-8pt]
&&\qquad\quad{} +\sst{i_1}{j_1}^3\sst{i_2}{j_2}^3\int_{t_1}^{t_2} ds_2\int
_0^{t_1} ds_1\,
\mu_0^2\bigl(\Ga^{(i_2j_2,0,i_1j_1,0,0)}_{
2;(s_1,t_1,s_2,t_3,t_4)}\psi\bigr)\nonumber\\
&&\qquad\quad{}  + \sst{i_1}{j_1}^3\sst{i_2}{j_2}^3\int_{t_1}^{t_2} ds_3\int
_0^{t_1} ds_2
\int_0^{s_2} ds_1\,
\mu_0\bigl(\Ga^{(12,i_2j_2,0,i_1j_1,0,0)}_{
1;(s_1,s_2,t_1,s_3,t_2,t_3,t_4)}\psi\bigr).\nonumber
\end{eqnarray}

Again from (\ref{eqtmsp}),
%
%
\begin{eqnarray}\label{eqfmp3}
&&\int_{t_2}^{t_3} ds\, \E\mu_{t_1}^3\bigl(\pi_1\Ga^{(0,12,0)}_{
2;(t_2-t_1,s-t_1,t_3-t_1,t_4-t_1)}\psi\bigr) \nonumber\\
&&\qquad=
\int_{t_2}^{t_3} ds\, \mu_0^3\bigl(\Ga
^{(0,0,12,0)}_{3;(t_1,t_2,s,t_3,t_4)}\psi\bigr)\nonumber\\[-8pt]\\[-8pt]
&&\qquad\quad{}+\sst{i}{j}^3\int_{t_2}^{t_3} ds_2\int_0^{t_1} ds_1
\,\mu_0^2\bigl(\Ga^{(ij,0,0,12,0)}_{2;(s_1,t_1,t_2,s_2,t_3,t_4)}\psi
\bigr)\nonumber\\
&&\qquad\quad{} + \sst{i}{j}^3\int_{t_2}^{t_3} ds_3\int_0^{t_1} ds_2 \int
_0^{s_2} ds_1\,
\mu_0\bigl(\Ga^{(12,ij,0,0,12,0)}_{1;(s_1,s_2,t_1,t_2,s_3,t_3,t_4)}\psi
\bigr)\nonumber.
\end{eqnarray}

From Lemma
\ref{lmmoment12sp},
%
%
\begin{eqnarray}\label{eqfmp4}
&& \sst{i}{j}^3\int_{t_1}^{t_2} ds_2 \int_{t_1}^{s_2} ds_1 \,\E\mu_{t_1}
^2\bigl(\pi_1\Ga
^{(12,ij,0,0)}_{1;(s_1-t_1,s_2-t_1,t_2-t_1,t_3-t_1,t_4-t_1)}\psi
\bigr)\nonumber\\
&&\qquad
=\sst{i}{j}^3\int_{t_1}^{t_2} ds_2 \int_{t_1}^{s_2} ds_1 \,\mu
_0^2\bigl(\Ga^{(0,12,
ij ,0,0)}_{2;(t_1,s_1,s_2,t_2,t_3,t_4)}\psi\bigr)\\
&&\qquad\quad{} +
\sst{i}{j}^3\int_{t_1}^{t_2} ds_3 \int_{t_1}^{s_3} ds_2 \int
_0^{t_1} ds_1\,
\mu_0\bigl(\Ga^{(12,0,12,
ij ,0,0)}_{1;(s_1,t_1,s_2,t_2,s_3,t_3,t_4)}\psi\bigr)\nonumber
\end{eqnarray}
and
%
%
\begin{eqnarray}\label{eqfmp5}
&&\int_{t_2}^{t_3} ds_2 \int_{t_1}^{t_2} ds_1 \,\E\mu_{t_1}^2\bigl(\pi
_1\Ga^{(12,
0,12,0)}_{1;(s_1-t_1,t_2-t_1,s_2-t_1,t_3-t_1,t_4-t_1)}\psi\bigr)\nonumber\\
&&\qquad= \int_{t_2}^{t_3} ds_2 \int_{t_1}^{t_2}
ds_1\,\mu_0^2\bigl(\Ga^ { (0 , 12 ,
0,12,0)}_{2;(t_1,s_1,t_2,s_2,t_3,t_4)}\psi\bigr) \\
&&\qquad\quad{} + \int_{t_2}^{t_3}
ds_3 \int_{t_1}^{t_2} ds_2 \int_0^{t_1}
ds_1\,\mu_0\bigl(\Ga^ { (12,0,12,
0,12,0)}_{1;(s_1,t_1,s_2,t_2,s_3,t_3,t_4)}\psi\bigr).\nonumber
\end{eqnarray}

Combining (\ref{eqfmp1}), (\ref{eqfmp2}), (\ref{eqfmp3}), (\ref{eqfmp4})
and (\ref{eqfmp5}), the desired formula follows.
\end{pf*}

\section{\texorpdfstring{Proof of Lemma \lowercase{\protect\ref{momentbound}}}{Proof of Lemma 3.10}}\label{appC}

We begin with some needed corollaries (of Theorem \ref{thmmomentsp}) and
lemmata.
\begin{coro}
\label{corbound1}For $i,j=1,2,3,4$, let $\phi^i_j\in
C^\infty_K(\R^{d})$ and define $\phi_i\in C^\infty_K(\R^{i\times
d})$ by
\[
\phi_i=\phi^i_1\otimes\cdots\otimes\phi^i_i.
\]
Then if $0\leq t_1\leq t_2\leq
t_3\leq t_4\leq T<\infty$,
\begin{eqnarray*}
\E\langle\phi_1,\mu_{t_1}\rangle&\leq& C(T)\Vert\phi_1\Vert
_\infty,
\\
\E\langle\phi_2,\mu_{t_1}\mu_{t_2}\rangle&\leq& C(T)\Vert\phi
_2\Vert_\infty,
\\
\E\langle\phi_3,\mu_{t_1}\mu_{t_2}\mu_{t_3}\rangle&\leq& C(T)\Vert
\phi_3\Vert_\infty
\end{eqnarray*}
and
\[
\E\langle\phi_4,\mu_{t_1}\mu_{t_2}\mu_{t_3}\mu_{t_4}\rangle\leq
C(T)\Vert\phi_4\Vert_\infty.
\]
\end{coro}
\begin{pf}
Since $\int dy q^{k}_t(x,y)=1$ for any $k\in\N$ and all $x\in\R
^{k\times
d}$, and since~$\mu_0$ is a finite measure having compact support,
this follows
immediately from Theorem \ref{thmmomentsp}.
\end{pf}
\begin{coro}
Equations (\ref{eqmixtmsp}) and (\ref{eqmixfmsp}) continue to hold
for $\phi\in S_{3\times d}$ and $\psi\in S_{4\times d}$.
\end{coro}
\begin{pf}
From Lemma \ref{lmsinglevariableapprox} there exist
$\{\phi_n\triangleq\sum_{k=1}^n\phi^1_{k}\otimes\phi^2_k\otimes
\phi^3_k\dvtx k\in\N\}
$ and
$\{\psi_n\triangleq\sum_{k=1}^n\psi^1_{k}
\otimes\psi^2_k\otimes\psi^3_k\otimes\psi^4_k\dvtx k\in\N\}$ such
that\vspace*{1pt}
$\phi^j_k,\psi^i_m\in C^\infty_K(\R^d)$, $i=1,2,3,4$, $j=1,2,3$,
$k,m\in\N$, and
$\lim_{n\rightarrow\infty}\phi_n=\phi$, $\lim_{n\rightarrow
\infty}\psi_n=\psi$, where the
convergence is uniform. For any $n,m\in N$, from equation
(\ref{eqmixfmsp}) and Corollary~\ref{corbound1}, it follows that
\[
\E\langle\vert\psi_n-\psi_m\vert,\mu_{t_1}\mu_{t_2}\mu
_{t_3}\mu_{t_4}\rangle\leq
C(T)\Vert\psi_n-\psi_m\Vert_\infty.
\]
Thus $\langle\psi_n,\mu_{t_1}\mu_{t_2}\mu_{t_3}\mu_{t_1}\rangle$ is
Cauchy in the
complete space $L^1(\P)$, and hence convergent. Uniform convergence of
$\psi_n$
implies
\[
\lim_{n\rightarrow\infty}\langle\psi_n,\mu_{t_1}\mu_{t_2}\mu
_{t_3}\mu_{t_4}\rangle=\langle\psi, \mu_{t_1}\mu_{t_2}\mu
_{t_3}\mu_{t_4}\rangle,\qquad\mbox{a.s.}
\]
Since the $L^1$ and a.s. limits must agree when they both exist,
\[
\lim_{n\rightarrow\infty}\E\langle\psi_n,\mu_{t_1}\mu_{t_2}\mu
_{t_3}\mu_{t_4}\rangle=\langle\psi,\mu_{t_1}\mu_{t_2}\mu
_{t_3}\mu_{t_4}\rangle.\vadjust{\goodbreak}
\]
Considering now the right-hand sides of equations (\ref{eqmixfmsp}) and
(\ref{eqmixtmsp}), by uniform
convergence, and since $\mu_0$ is finite with compact support, the desired
convergence is shown.
\end{pf}

For ease in reading, we introduce the notation
%
%
\begin{equation}
\label{fnXi}\quad
\Xi^{\mathbf{x}_{m - 1}}_{\ell;\mathbf{s}_{m}} \triangleq
Q^{\ell}_{s_1}\zeta(x_1)Q^{\ell(x_1)}_{s_2-s_1}\zeta(x_2)Q^{\ell
(x_2)}_{s_3-s_2}
\cdots
\zeta(x_{m-2})Q^{3}_{s_{m-1}-s_{m-2}}\pi_1Q^2_{s_m-s_{m-1}},\hspace*{-32pt}
\end{equation}
where $\mathbf{x}_m$, $\mathbf{s}_m$, $\ell$, $\ell(x)$ and $\zeta$
are defined as in
(\ref{fnGamma}), with the convention that $x_{m-1}=0$ and $\ell(m-1)=2$.
\begin{pf*}{Proof of Lemma \ref{momentbound}}
Throughout this proof, the norm on $L^p$ will be denoted by $\Vert
\cdot\Vert_p$.
From the moment equation (\ref{eqmixfmsp}) and the preceding
corollary, it
follows that
\[
\int_0^{t_3} dt_2 \int_0^{t_2} dt_1\, \E\langle\vp,\mu_{t_1}\mu
_{t_2}\mu_{t_3}^2 \rangle\leq C\sum_{k=1}^{14}J_k(t_1,t_2,t_3),
\]
where the definition of each $J_k$ is implicit in equation (\ref{eqmixfmsp}).

To begin, note that from Dynkin (\citeyear{dynkin65}), Theorem 0.5,
for $n\in\N$, $x,\break y\in\R^{n\times d}$,
\[
q^n_t(x,y)\leq C p^n_{\io t}(x,y),
\]
where $C$ and $\io$ are constants, and $p^n_\cdot=\prod_{i=1}^n
p_\cdot$, where $p_\cdot$ is the Brownian transition function on $\R
^{d}$. It thus follows that
%
%
\begin{eqnarray}\label{F1bound}
&&\Xi^{(0,0)}_{4;(t_1-s,t_2-s,t_3-s)}\vp(x)\nonumber\\
&&\qquad\leq
C\int da\, p_{\io(t_1-s)}(x_1,a_1)p_{\io(t_2-s)}(x_2,a_2)\\
&&\hspace*{17.5pt}\qquad\quad{}\times p_{\io(t_3-s)}(x_3,
a_3)p_{\io(t_3-s)}(x_4,a_4)\vp(a)\nonumber
\end{eqnarray}
for all $ x\in\R^{4\times d}, s\in[0,t_1]$, $a=(a_1,a_2,a_3,a_4)$.

Using inequality (\ref{F1bound}), it follows that
\begin{eqnarray*}
&&
\bigl\langle\Xi^{(0,0)}_{4;(t_1,t_2,t_3)}\vp,\mu_0^4\bigr\rangle\\
&&\qquad\leq
C \int\mu_0^4(dx) \int da\, p_{\io t_3}(x_4,a_4)\prod
_{i=1}^3 p_{\io
t_i}(x_i,a_i)\vp(a_1,a_2,a_3,a_4)\\
&&\qquad\leq C\int\mu_0(dx_3)\mu_0(dx_4) \int da_3\,da_4\, p_{\io
t_3}(x_3,a_3)p_{\io
t_3}(x_4,a_4) \int da_1\, \phi(a_1-a_3)\\
&&\qquad\quad{} \times\int da_2\,\phi(a_2-a_4) \int\mu_0(dx_1)p_{\io
t_1}(x_1,a_1)\int\mu_0(dx_2) p_{\io t_2}(x_2,a_2)\\
&&\qquad\leq C\Vert\phi\Vert_1^2,
\end{eqnarray*}
and thus, since $d\leq3$,
\[
\int_0^{t_3} dt_2 \int_0^{t_2} dt_1 \bigl\langle\Xi
^{(0,0)}_{4;(t_1,t_2,t_3)}\vp, \mu_0^4\bigr\rangle\leq C(T)\Vert\phi
\Vert_1^2.
\]
Let $\{i',j'\}=\{1,2,3,4\}\setminus\{i,j\}, i'< j'$, then again from
(\ref{F1bound}),
\begin{eqnarray*}
\hspace*{-5pt}&&\sst{i}{j}^4 \int_0^{t_1} ds \bigl\langle\Xi
^{(ij,0,0)}_{3;(s,t_1,t_2,t_3)}\vp, \mu^3_0 \bigr\rangle\\
\hspace*{-5pt}&&\qquad \leq C\sstl{i}{j}^4 \int_0^{t_1} ds \int\mu_0^3(dx) \int
dy\, p_{\io
s}(x_1,y) \int da_1\,da_2\,da_3\,da_4\, p_{\io(t_i-s)}(y,a_i)\\
\hspace*{-5pt}&&\hspace*{28.8pt}\qquad\quad{} \times p_{\io(t_{j\mn3}-s)}(y,a_j) p_{\io t_{i'}}(x_2,a_{i'})p_{\io
t_{j'\mn3}}(x_3,a_{j'})\phi(a_1-a_3)\phi(a_2-a_4)
\end{eqnarray*}
and so, with a some applications of the Kolmogorv--Chapman equation to
the above
expression,
\[
\sst{i}{j}^4 \int_0^{t_1} ds \bigl\langle\Xi
^{(ij,0,0)}_{3;(s,t_1,t_2,t_3)}\vp, \mu^3_0 \bigr\rangle\leq C(T)\Vert
\phi\Vert_1^2+C\int_0^{t_1} ds (t_3-s)^{-d/2}.
\]
Since $d\leq3$, it follows that
\[
\sst{i}{j}^4 \int_0^{t_3} dt_2 \int_0^{t_2} dt_1 \int_0^{t_1}
ds \bigl\langle\Xi^{(ij,0,0)}_{3;(s,t_1,t_2,t_3)}\vp,\mu^3_0\bigr\rangle
\leq C(T)\Vert\phi\Vert_1^2.
\]

This next case becomes quite a bit more complicated, so we explain with more
detail. Consider
\[
\sst{i}{j}^4\sst{n}{m}^3\int_0^{t_1} ds_2 \int_0^{s_2} ds_1
\bigl\langle\Xi^{ (nm,ij,0,0)}_{2;(s_1,s_2,t_1,t_2,t_3)}\vp,\mu
_0^2\bigr\rangle,
\]
wherein the presence of both $\Phi_{nm}$ and $\Phi_{ij}$ greatly
increase the
number of cases. In bounding, we again may assume, with the addition of a
multiplicative constant to the bound, that $i<j$. Note first that when
$m+n\neq
6-i$ it will follow that either
\begin{eqnarray*}
&&\Xi^{(nm,ij,0,0)}_{2;(s_1,s_2,t_1,t_2,t_3)}\vp(x_1,x_2)\\
&&\qquad\leq C\int dy\, p_{\io
s_1}(x_1,y) \int dz\, p_{\io(s_2-s_1)}(y,z)\\
&&\qquad\quad{}\times\int dw_1\,dw_2\,dw_3\,dw_4\, p_{
\io(t_i-s_2)}(z,w_i)\\
&&\hspace*{20.6pt}\qquad\quad{}\times
p_{\io(t_{j\mn3}-s_2)}(z,w_j)p_{
\io(t_{i'}-s_1)}(y,w_{i'})p_{\io t_{j'\mn3}}(x_2,w_{j'})\vp(w)
\end{eqnarray*}
or
\begin{eqnarray*}
&&\Xi^{
(nm,ij,0,0)}_{2;(s_1,s_2,t_1,t_2,t_3)}\vp(x_1,x_2)\\[-2pt]
&&\qquad\leq C\int dy\, p_{\io
s_1}(x_1,y) \int dz\, p_{\io(s_2-s_1)}(y,z)\\[-2pt]
&&\qquad\quad{}\times\int dw_1\,dw_2\,dw_3\,dw_4\, p_{
\io(t_i-s_2)}(z,w_i)\\[-2pt]
&&\hspace*{20.6pt}\qquad\quad{}\times
p_{\io(t_{j\mn3}-s_2)}(z,w_j)p_{
\io(t_{j'\mn3}-s_1)}(y,w_{j'})p_{\io t_{i'}}(x_2,w_{i'})\vp(w),
\end{eqnarray*}
where again $\{i',j'\}=\{1,2,3,4\}\setminus\{i,j\}$, with $i'<j'$. In
the case
that $m+n=6-i$, we have the bound
\begin{eqnarray*}
&&\Xi^{
(nm,ij,0,0)}_{2;(s_1,s_2,t_1,t_2,t_3)}\vp(x_1,x_2)\\[-2pt]
&&\qquad\leq C\int dy\, p_{\io s_1}(x_1,y) \int dz\, p_{\io
s_2}(x_2,z)\\[-2pt]
&&\qquad\quad{}\times\int dw_1\,dw_2\,dw_3\,dw_4\, p_{\io(t_i-s_2)}(z,w_i)\\[-2pt]
&&\hspace*{20.6pt}\qquad\quad{}\times
p_{\io(t_{j\mn3}-s_2)}(z,w_j)p_{
\io(t_{i'}-s_1)}(y,w_{i'})p_{\io(t_{j'\mn3}-s_1)}(y,w_{j'})\vp(w).
\end{eqnarray*}
It thus follows that
\[
\sst{i}{j}^4\sst{n}{m}^3\int_0^{t_1} ds_2 \int_0^{s_2} ds_1
\bigl\langle\Xi^{ (nm,ij,0,0)}_{2;(s_1,s_2,t_1,t_2,t_3)}\vp,\mu
_0^2\bigr\rangle\leq C(A_1+A_2+A_3),
\]
where
\begin{eqnarray*}
A_1&=&\mathop{\mathop{\sum_{i,j,i',j'=1}}_{i'\neq i,j'\neq
j}}_{i<j,i'<j'}^4\mathop{\mathop{\sum_{n,m=1}}_{n<m}}_{6-n-m\neq
i}^3\int_0^{t_1} ds_2 \int_0^{s_2} ds_1 \int\mu_0(dx_1) \int
dy \,p_{\io
s_1}(x_1,y)\\[-4pt]
&&\hspace*{75.6pt}{}
\times\int dz\, p_{\io(s_2-s_1)}(y,z)\\[-2pt]
&&\hspace*{75.6pt}{}
\times\int dw_i\,dw_j\,dw_{i'}\, p_{\io
(t_{i'}-s_1)}
(y,w_{i'})p_{\io(t_i-s_2)}(z,
w_i)\\[-2pt]
&&\hspace*{36.2pt}\hspace*{61pt}{}\times
p_{\io(t_{j\mn3}-s_2)}(z,w_j)\int dw_{j'}\,
\phi(w_1-w_3)\phi(w_2-w_4)\\[-2pt]
&&\hspace*{75.6pt}{}\times\int\mu_0(dx_2) p_{\io
t_{j'\mn3}}(x_2,w_{j'}),
\\[-2pt]
A_2&=&\mathop{\mathop{\sum_{i,j,i',j'=1}}_{i'\neq i,j'\neq
j}}_{i<j,i'<j'}^4\mathop{\mathop{\sum_{n,m=1}}_{n<m}}_{6-n-m\neq
i}^3\int_0^{t_1} ds_2 \int_0^{s_2} ds_1 \int\mu_0(dx_1) \int
dy\, p_{\io
s_1}(x_1,y)\\[-3pt]
&&\hspace*{75.6pt}{}
\times\int dz\, p_{\io(s_2-s_1)}(y,z)\\[-0.8pt]
&&\hspace*{75.6pt}{}\times\int dw_i\,dw_j\,dw_{j'}\, p_{\io
(t_{j'\mn3}
-s_1)}(y,w_{j'})p_{\io(t_i-s_2)
}(z,w_i)\\[-0.8pt]
&&\hspace*{75.6pt}\hspace*{21.22pt}{}\times
p_{\io(t_{j\mn3}-s_2)}(z,w_j)\int dw_{i'}\,
\phi(w_1-w_3)\phi(w_2-w_4)\\[-0.8pt]
&&\hspace*{75.6pt}{}\times\int\mu_0(dx_2) p_{\io
t_{i'}}(x_2,w_{i'})
\end{eqnarray*}
and
\begin{eqnarray*}
A_3&=&\mathop{\mathop{\sum_{i,j,i',j'=1}}_{i'\neq i,j'\neq
j}}_{i<j,i'<j'}^4\mathop{\mathop{\sum_{n,m=1}}_{n<m}}_{6-n-m=
i}^3\int_0^{t_1} ds_2 \int_0^{s_2}
ds_1 \int\mu_0(dx_1)\mu_0(dx_2) \\[-3pt]
&&\hspace*{75.6pt}{}\times\int dy \,p_{\io
s_1}(x_1,y)\int dz\, p_{\io
s_2}(x_2,z)\\[-0.8pt]
&&\hspace*{75.6pt}{}\times\int dw_1\,dw_2\,dw_3\,dw_4\, p_{\io
(t_{i'}-s_1)}(y,w_{i'})\\[-0.8pt]
&&\hspace*{75.6pt}\hspace*{21.6pt}{}\times p_{\io(t_{j'\mn3
}
-s_1)}(y,w_{j'})
p_{\io(t_i-s_2)}(z,w_i)\\[-0.8pt]
&&\hspace*{75.6pt}\hspace*{21.6pt}{}\times p_{\io(t_{j\mn3}-s_2)}(z,w_j)\vp(w_1,w_2,w_3,w_4).
\end{eqnarray*}

For the first of the above three terms, the process is as follows. Bound
$\mu_0(dx_2)$ by $\Vert m\Vert_\infty \,dx_2$, integrate $p_{\io
t_{j'\mn
3}}(x_2,w_j')$ with
respect to $dx_2$, then $\phi$ with respect to $dw_{j'}$. In doing so,
we may
then integrate out one of the remaining $w_i,w_j$ or $w_j'$. In what remains,
if $(i,j)\neq(1,2)$ there will be the term $p_{\io
(t_3-s_2)}(z,w_\cdot)$, or if
$(i,j)=(1,2)$, the term $p_{\io(t_3-s_1)}(y,w_3)$. In either case,
bound the
respective term by $C(t_3-s_.)^{-d/2}$. This allows for the integration
of the
second $\phi$.

For the second of the two above terms,
bound $\mu_0(dx_2)$ by $\Vert m\Vert_\infty\, dx_2$, integrate $p_{\io
t_{i'\mn
3}}(x_2,w_i')$
with respect to $dx_2$, then $\phi$ with respect to $dw_{i'}$. In
doing so, we
may then integrate out one of the remaining $w_i,w_j$ or $w_j'$. In what
remains, if $(i,j)\in\{(1,2),(1,4),(2,3)\}$ there will be the term
$p_{\io(t_3-s_1)}(y,w_{j'})$, otherwise there will exist the term
$p_{\io(t_3-s_2)}(z,w_{\cdot})$. In either case, bound the respective
term by
$C(t_3-s_.)^{-d/2}$. This allows for the integration of the second
$\phi$.

For the third and final term,
if $(i,j)\notin\{(1,2),(3,4)\}$, there will exist the terms
$p_{\io(t_3-s_2)}(z,w_j)$ and $p_{\io(t_3-s_1)}$, which are bounded,
respectively,
by $C(t_3-s_2)^{-d/2}$ and $C(t_3-s_1)^{-d/2}$. When\vspace*{1pt} $(i,j)=(1,2)$ we bound
$p_{\io(t_3-s_1)}(y,w_3)$ and $p_{\io(t_2-s_2)}(z,w_2)$, respectively, by
$C(t_3-s_1)^{-d/2}$ and $C(t_2-s_2)^{-d/2}$. Finally, when
$(i,j)=(3,4)$, bound
the terms $p_{\io(t_3-s_2)}(z,w_3)$ and
$p_{\io(t_2-s_1)}(y,w_2)$,\vspace*{1pt}
respectively,
by $C(t_3-s_2)^{-d/2}$ and $C(t_2-s_1)^{-d/2}$. This allows for the desired
integration of $\phi(w_1-w_3)\phi(w_2-w_4)$.

Combining the above, and since $d\leq3$, we arrive at the
bound\vspace*{-2.5pt}
\[
\sst{i}{j}^4\sst{n}{m}^3\int_0^{t_3} dt_2 \int_0^{t_2}
dt_1 \int_0^ { t_1 }
ds_1 \int_0^{s_1} ds_2\bigl\langle\Xi^{
(nm,ij,0,0)}_{2;(s_1,s_2,t_1,t_2,t_3)}\vp,\mu_0^2\bigr\rangle\leq
C(T)\Vert\phi\Vert_1^2.
\]

Considering the next case, note first the similarities in the respective
corresponding particle pictures of this and the previous case. This
case can be
seen as a modification of the previous case in which the two original particles
were both born from a common ancestor. Thus, arguing as before, we
arrive at the
bound\vspace*{-2.5pt}
\[
\sst{i}{j}^4\sst{n}{m}^3 \int_0^{t_1} ds_3 \int_0^{s_3}
ds_2 \int_0^ { s_2 }
ds_1\bigl\langle\Xi^{(12,nm,ij,0,0)}_{1;(s_1,s_2,s_3,t_1,t_2,t_3)}\vp
,\mu_0 \bigr\rangle\leq C (B_1+B_2+B_3),
\]
where for the $B_k$ we have
\begin{eqnarray*}
B_1&=&\mathop{\mathop{\sum_{i,j,i',j'=1}}_{i'\neq i,j'\neq
j}}_{i<j,i'<j'}^4\mathop{\mathop{\sum_{n,m=1}}_{n<m}}_{6-n-m\neq
i}^3 \int_0^{t_1} ds_3 \int_0^{s_3} ds_2 \int_0^{s_2}
ds_1 \int\mu_0(dx) \int dy\, p_{\io
s_1}(x,y)\\[-4.5pt]
&&\hspace*{75.6pt}{}
\times\int dz \,p_{\io(s_2-s_1)}(y,z)\int dw \,p_{\io(s_3-s_2)}(z,
w)\\[-2pt]
&&\hspace*{75.6pt}{}\times\int dv_1\,dv_2\,dv_3\,dv_4\,
p_{\io(t_{j'\mn3}-s_1)}(y,v_{j'})p_{\io(t_{i'}-s_2)}(z,v_{i'})\\
&&\hspace*{75.6pt}\hspace*{21.9pt}{}\times
p_{\io(t_i-s_3)}(w,v_i)p_{\io(t_{j\mn3}-s_3)}(w,
v_j)\\[-1pt]
&&\hspace*{75.6pt}\hspace*{21.9pt}{}\times\phi(v_1-v_3)\phi(v_2-v_4),
\\
B_2&=&\mathop{\mathop{\sum_{i,j,i',j'=1}}_{i'\neq i,j'\neq
j}}_{i<j,i'<j'}^4\mathop{\mathop{\sum_{n,m=1}}_{n<m}}_{6-n-m\neq
i}^3 \int_0^{t_1} ds_3 \int_0^{s_3} ds_2 \int_0^{s_2}
ds_1 \int\mu_0(dx) \int dy\, p_{\io
s_1}(x,y)\\
&&\hspace*{75.6pt}{}\times \int dz \,p_{\io(s_2-s_1)}(y,z)\int dw \,
p_{\io(s_3-s_2)}(z,w)\\
&&\hspace*{75.6pt}{}\times\int dv_1\,dv_2\,dv_3\,dv_4\,
p_{\io(t_{i'}-s_1)}(y,v_{i'})p_{\io(t_{j'\mn3}-s_2)}(z,v_{j'})\\
&&\hspace*{75.6pt}\hspace*{21.9pt}{}\times
p_{\io(t_i-s_3)}(w,v_i)p_{\io(t_{j\mn3}-s_3)}(w,
v_j)\\
&&\hspace*{75.6pt}\hspace*{21.9pt}{}\times\phi(v_1-v_3)\phi(v_2-v_4)
\end{eqnarray*}
and
\begin{eqnarray*}
B_3&=&\mathop{\mathop{\sum_{i,j,i',j'=1}}_{i'\neq i,j'\neq
j}}_{i<j,i'<j'}^4\mathop{\mathop{\sum_{n,m=1}}_{n<m}}_{6-n-m=
i}^3 \int_0^{t_1} ds_3 \int_0^{s_3} ds_2 \int_0^{s_2}
ds_1 \int\mu_0(dx)\int dy\, p_{\io s_1}(x,y)\\
&&\hspace*{75.6pt}{}\times
\int dz_1\,dz_2\, p_{\io(s_2-s_1)}(y,z_1)p_{\io(s_3-s_1)}(y,
z_2)\\
&&\hspace*{75.6pt}{}\times \int dv_1\,dv_2\,dv_3\,dv_4\, p_{\io(t_i-s_3)}(z_2,v_i)
p_{\io(t_{j\mn3}-s_3)}(z_2,v_j)\\
&&\hspace*{75.6pt}\hspace*{21.9pt}{}\times p_{\io(t_{i'}-s_2)}(z_1,v_{i'})p_{\io
(t_{j'\mn3}
-s_2)}(z_1,v_{j'})\\
&&\hspace*{75.6pt}\hspace*{21.9pt}{}\times\phi(v_1-v_3)\phi(v_2-v_4).
\end{eqnarray*}
Thus,
\begin{eqnarray*}
&&B_1+B_2+B_3\\
&&\qquad\leq
C\Vert\phi\Vert_1^2 \int_0^{t_1} ds_3 \int_0^{s_3}
ds_2 \\
&&\qquad\quad{}\times \int_0^{
s_2 } ds_1
[(t_3-s_2)^{-d/2}(t_2-s_2)^{-d/2}+ (t_3-s_2)^{-{d/2}
} (t_3-s_3)^{-{d/2}
}\\
&&\qquad\quad\hspace*{49pt}{} + (t_2-s_2)^{-d/2} (t_3-s_3)^{-d/2} +
(t_3-s_2)^{-{d
}/ { 2 } } (t_2-s_3)^ { -{d/2}}],
\end{eqnarray*}
where the above bounds are obtained similarly to the previous case.
And so, since $d\leq3$,
\begin{eqnarray*}
&&\sst{i}{j}^4\sst{n}{m}^3\int_0^{t_3} dt_2 \int_0^{t_2} dt_1
\int_0^{t_1}
ds_3 \int_0^{s_3} ds_2\int_0^{s_2} ds_1\bigl\langle\Xi
^{(12,nm,ij,0,0)}_{1;(s_3, s_2,s_1,t_1,t_2,t_3)}\vp, \mu_0 \bigr\rangle
\nonumber\\
&&\qquad \leq C(T)\Vert\phi\Vert_1^2.
\end{eqnarray*}

This takes care of four of the fourteen $J_k$, and we consider now the
next three
integrals which are dependent upon the expression
%
%
\begin{eqnarray}\label{F2bound}
&&\Xi^{(0,ij,0)}_{3;(t_1-s_1,s_2-s_1,t_2-s_1,t_3-s_1)}
\vp(x) \nonumber\\[-1.2pt]
&&\qquad\leq
C\int db\, p_{\io(s_2-s_1)}(x_3,b) \nonumber\\[-1.2pt]
&&\hspace*{17.8pt}\qquad\quad{}\times\int da_1\,da_2\,da_3\,da_4\, p_{\io
(t_1-s_1)}(x_1
,a_1)p_{\io(t_{4-i}-s_1)}(x_2,a_{7-i-j})\\[-1.2pt]
&&\qquad\quad\hspace*{40.2pt}{}\times
p_{\io(t_{i+1}-s_2)}(b,a_{i+1})p_{\io(t_3-s_2)}(b,a_{j+1})\nonumber\\[-1.2pt]
&&\qquad\quad\hspace*{40.2pt}{}\times\phi(a_1-a_3)\phi(a_2-a_4)\nonumber
\end{eqnarray}
for all $x\in\R^{3\times d}$, $0\leq s_1\leq t_1$, $t_1\leq s_2\leq
t_2$ and
$i,j=1,2,3$, $i\neq j$. In the above $z^{ij}_\cdot$ refers to the particular
arrangement of $z_1,z_2$ given the pair $(i,j)$. Applying (\ref
{F2bound}) now
gives
\begin{eqnarray*}
&&\sst{i}{j}^3\int_{t_1}^{t_2} ds \bigl\langle\Xi
^{(0,ij,0)}_{3;(t_1,s,t_2,t_3)}\vp, \mu_0^3 \bigr\rangle\\[-2.5pt]
&&\qquad\leq C
\sstl{i}{j}^3\int_{t_1}^{t_2} ds \int\mu_0(dx_2)\mu_0(dx_3)
\int dy\, p_{
\io s}(x_3,y)\\[-2.5pt]
&&\qquad\quad\hspace*{28.6pt}{}\times\int dz_{7-i-j}\,dz_{i+1} p_{\io t_{4-i}}(x_2,z_{7-i-j})
p_{\io(t_{i+1}-s)}(y,z_{i+1})\\[-1.2pt]
&&\qquad\quad\hspace*{28.6pt}{}\times\int dz_{j+1}\, p_{\io
(t_3-s)}(y,z_{j+1})\int dz_1\, \phi(z_1-z_3)\phi(z_2-z_4)\\[-1.2pt]
&&\qquad\quad\hspace*{28.6pt}{}\times\int\mu_0(dx_1)
p_{\io
t_1}(x_1,z_1)
\end{eqnarray*}
and, by integrating out terms and using known bounds,
\[
\sst{i}{j}^3\int_{t_1}^{t_2} ds \bigl\langle\Xi
^{(0,ij,0)}_{3;(t_1,s,t_2,t_3)}\vp, \mu_0^3 \bigr\rangle\leq C\Vert\phi
\Vert_1^2\int_{t_1}^{t_2} ds
[1+(t_3-s)^{-d/2}].
\]
And so, since $d\leq3$,
\[
\sst{i}{j}^3\int_0^{t_3} dt_2 \int_0^{t_2} dt_1 \int
_{t_1}^{t_2} ds\bigl\langle\Xi^{(0,ij,0)}_{3;(t_1,s,t_2,t_3)}\vp,\mu
_0^3\bigr\rangle\leq C(T)\Vert\phi\Vert_1^2.
\]

Again from (\ref{F2bound}),
\[
\sst{i}{j}^3\sst{n}{m}^3\int_{t_1}^{t_2} ds_2 \int_0^{t_1} ds_1
\bigl\langle\Xi^{(nm,0,ij,0)}_{2;(s_1,t_1,s_2,t_2,t_3)}\vp, \mu
_0^2\bigr\rangle\leq C(C_1+C_2+C_3),
\]
where for the $C_k$ we have
\begin{eqnarray*}
C_1&=&\sstl{i}{j}^3\int_{t_1}^{t_2} ds_2 \int_0^{t_1}
ds_1 \int\mu_0(dx_1) \int dy\, p_{\io
s_1}(x_1,y)\int dz\, p_{\io(s_2-s_1)}(y,z)\\[-3pt]
&&\hspace*{18.8pt}{} \times\int dw_1\,dw_{i+1}\,dw_{j+1}\,
p_{\io(t_1-s_1)}(y,w_1)\\[-3pt]
&&\qquad\quad\hspace*{6.5pt}{}\times
p_{\io(t_{i+1}-s_2)}(z,w_{i+1})p_{\io(t_3-s_2)}(z,w_{j+1})\\[-3pt]
&&\hspace*{18.8pt}{}\times\int dw_{7-i-j}\,
\phi(w_1-w_3)\phi(w_2-w_4)\\[-3pt]
&&\hspace*{18.8pt}{} \times\int\mu_0(dx_2) p_{\io(t_{(7-i-j)\mn
3})}(x_2,w_{
7-i-j})\\[-3pt]
&\leq&
C\Vert\phi\Vert_1\int_{t_1}^{t_2} ds_2 \int_0^{t_1} ds_1 (t_3-s_2)^{-d/2},
\\[-3pt]
C_2&=&\sstl{i}{j}^3\int_{t_1}^{t_2} ds_2 \int_0^{t_1}
ds_1 \int\mu_0(dx_1) \int dy\, p_{\io
s_1}(x_1,y)\int dz\, p_{\io(s_2-s_1)}(y,z)\\[-3pt]
&&\hspace*{18.8pt}{} \times\int dw_2\,dw_3\,dw_4\, p_{\io(t_{(7-i-j)\mn3}-s_1)}(y,w_{7-i-j})
p_{\io(t_{i+1}-s_2)}(z,w_{i+1})\\[-3pt]
&&\qquad\quad\hspace*{6.5pt}{}\times
p_{\io(t_3-s_2)}(z,w_{j+1})\phi(w_2-w_4)\int dw_1\,
\phi(w_1-w_3)\\[-3pt]
&&\hspace*{18.8pt}{} \times\int\mu_0(dx_2) p_{\io t_1}(x_2,w_1)\\[-3pt]
&\leq&
C\Vert\phi\Vert_1\int_{t_1}^{t_2} ds_2 \int_0^{t_1} ds_1
[(t_3-s_1)^{-d/2}
+(t_3-s_2)^{-d/2}]
\end{eqnarray*}
and
\begin{eqnarray*}
C_3&=&C\sstl{i}{j}^3
\int_{t_1}^{t_2} ds_2 \int_0^{t_1}
ds_1 \int\mu_0(dx_1)\mu_0(dx_2)\\[-3pt]
&&\hspace*{29.1pt}{} \times \int dy\, p_{\io
s_1}(x_1,y)\int dz\, p_{\io s_2}(x_2,z)\\[-3pt]
&&\hspace*{29.1pt}{} \times\int dw_1\,dw_2\,dw_3\,dw_4\, p_{\io(t_1-s_1)}(y,w_1)
p_{\io(t_{(7-i-j)\mn3}-s_1)}(y,w_{7-i-j})\\[-3pt]
&&\hspace*{49.3pt}{} \times p_{\io(t_{i+1}-s_2)}(z,w_{i+1})
p_{\io(t_3-s_2)}(z,w_{j+1})\phi(w_1-w_3)\phi(w_2-w_4)\\[-3pt]
&\leq&
C\Vert\phi\Vert_1\int_{t_1}^{t_2} ds_2 \int_0^{t_1} ds_1
[(t_3-s_1)^{-d/2}
(t_3-s_2)^{-d/2}+(t_3-s_2)^{-d/2}].
\end{eqnarray*}
Therefore,
\begin{eqnarray*}
&& C_1+C_2+C_3\\[-3pt]
&&\qquad\leq C\Vert\phi\Vert_1^2
\int_{t_1}^{t_2} ds_2 \int_0^{t_1} ds_1
[(t_3-s_1)^{-{d}/{2
}}
+ (t_3-s_1)^{-{d}/{2
}}(t_3-s_2)^{-{d}/{2
}} \\[-3pt]
&&\qquad\quad\hspace*{220pt}{}+ (t_3-s_2)^{-{d}/{2
}}].
\end{eqnarray*}
Since $d\leq3$,
\[
\sst{i}{j}^3\sst{n}{m}^3\int_0^{t_3} dt_2 \int_0^{t_2}
dt_1 \int_{
t_1}^{t_2} ds_2 \int_0^{t_1} ds_1 \bigl\langle\Xi^{(nm,0,ij,0)}_{2;(s_1,
t_1,s_2,t_2,t_3)}\vp, \mu_0^2 \bigr\rangle\leq C(T)\Vert\phi\Vert_1^2.
\]
After one final application of (\ref{F2bound}),
\[
\sst{i}{j}^3 \sst{n}{m}^3 \int_{t_1}^{t_2} ds_3 \int_0^{t_1}
ds_2 \int_0^{s_2} ds_1 \bigl\langle\Xi^{(12,nm,0,ij,0)}_{1;(s_1,s_2,
t_1,s_3,t_2,t_3)}\vp, \mu_0 \bigr\rangle\leq C(D_1+D_2+D_3),
\]
where for the $D_k$ we have
\begin{eqnarray*}
D_1&=&
C\sst{i}{j}^3 \int_{t_1}^{t_2} ds_3 \int_0^{t_1} ds_2
\int_0^{
s_2 }
ds_1 \\
&&\hspace*{28pt}{}\times\int\mu_0(dx) \int dy\, p_{\io s_1}(x,y)
\int dz\, p_{\io(s_2-s_1)}(y,z)
\int dw\, p_{\io(s_3-s_2)}(z,w) \\
&&\hspace*{28pt}{}\times\int dv_1\,dv_2\,dv_3\,dv_4\, p_{
\io(t_1-s_2) } (z , v_1)
p_{\io(t_{(7-i-j)\mn3}-s_1)}(y,v_{7-i-j})\\
&&\qquad\quad\hspace*{16.4pt}{}\times
p_{\io(t_{i+1}-s_3)}(w,v_{i+1})p_{\io(t_3-s_3)}(w,v_{j+1}
)\phi(v_1-v_3)\phi(v_2-v_4)\\
&\leq&
C \Vert\phi\Vert_1^2 \int_{t_1}^{t_2} ds_3 \int_0^{t_1}
ds_2 \int_0^ { s_2 } ds_1
[ (t_3-s_1)^{-d/2}(t_3-s_3)^{-d/2} \\
&&\qquad\quad\hspace*{105pt}{}+ (t_2-s_1)^{
-{d}/{2}} (t_3-s_3)^ { -{d}/{2} } ],
\\
D_2&=&
C \sst{i}{j}^3 \int_{t_1}^{t_2} ds_3 \int_0^{t_1} ds_2 \int_0^{s_2}
ds_1 \int\mu_0(dx)\int dy\, p_{\io s_1}(x,y)\\
&&\hspace*{28.8pt}{} \times
\int dz\, p_{\io(s_2-s_1)}(y,z)
\int dw\, p_{\io(s_3-s_2)}(z,w)\\
&&\hspace*{28.8pt}{} \times\int dv_1\,dv_2\,dv_3\,dv_4\, p_{\io
(t_{(7-i-j)\mn3
}
-s_2)}(z,v_{7-i-j}) p_{\io(t_1-s_1)}(y,v_1)\\
&&\hspace*{48.8pt}{} \times
p_{\io(t_{i+1}-s_3)}(w,v_{i+1})p_{\io(t_3-s_3)}(w,v_{j+1}
)\phi(v_1-v_3)\phi(v_2-v_4)\\
&\leq&
C \Vert\phi\Vert_1^2 \int_{t_1}^{t_2} ds_3 \int_0^{t_1}
ds_2 \int_0^ { s_2
} ds_1 [ (t_3-s_2)^{-d/2}(t_3-s_3)^{-d/2}\\
&&\hspace*{138.4pt}{} + (t_2-s_2)^ { -{d}/{2} } (t_3-s_3)^ { -{d}/{2} }
]
\end{eqnarray*}
and
\begin{eqnarray*}
D_3&=&
C \sst{i}{j}^3 \int_{t_1}^{t_2} ds_3 \int_0^{t_1}
ds_2 \int_0^ { s_2 }
ds_1 \int\mu_0(dx)\int dy\, p_{\io s_1}(x,y)\\
&&\hspace*{29.1pt}{}\times
\int dz\, p_{\io(s_2-s_1)}(y,z)\int dw\,
p_{\io(s_3-s_1)}(y,w) \\
&&\hspace*{29.1pt}{}\times\int dw_1\,dw_2\,dw_3\,dw_4\, p_{\io(t_1-s_2)}(z,v_1)
p_{\io(t_{(7-i-j)\mn3}-s_2)}(z,v_{7-i-j})\\
&&\hspace*{50pt}{}\times
p_{\io(t_{i+1}-s_3)}(w,v_{i+1})p_{\io(t_3-s_3)}(w,v_{j+1}
)\phi(v_1-v_3)\phi(v_2-v_4)\\
&\leq&
C \Vert\phi\Vert_1^2 \int_{t_1}^{t_2} ds_3 \int_0^{t_1}
ds_2 \int_0^ { s_2 }
ds_1 [ (t_3-s_2)^{-d/2}(t_3-s_3)^{-
{d}/{2}}\\
&&\hspace*{138.1pt}{}+(t_2-s_2)^{
-{d}/{2} } (t_3-s_3)^ { -{d}/{2}} ] .
\end{eqnarray*}
Thus,
\begin{eqnarray*}
&&D_1+D_2+D_3\\
&&\qquad
\leq C \Vert\phi\Vert_1^2\sum_{k=1}^2 \int_{t_1}^{t_2}
ds_3 \int_0^{
t_1 } ds_2 \\
&&\qquad\quad\hspace*{45.2pt}{}\times\int_0^ {
s_2} ds_1\bigl[(t_3-s_3)^{-d/2}\bigl((t_3-s_k)^{-d/2}
+(t_2-s_k)^ { -{d}/{2} }
\bigr)\bigr].
\end{eqnarray*}
Therefore, since $d\leq3$,
\begin{eqnarray*}
&&\sst{i}{j}^3 \sst{n}{m}^3 \int_0^{t_3} dt_2 \int_0^{t_2}
dt_1 \int_ { t_1 } ^ { t_2
} ds_3 \int_0^{t_1} ds_2 \int_0^{s_2} ds_1 \bigl\langle\Xi^{(12, nm , 0
, ij , 0) } _ { 1;(s_1,s_2,t_1,s_3,t_2,t_3)} , \mu_0 \bigr\rangle\\
&&\qquad\leq
C(T) \Vert\phi\Vert_1^2.
\end{eqnarray*}
Thus seven of the fourteen $J_k$ are now shown to have the desired
bound, we
continue with three more of the $J_k$.
%
%
\begin{eqnarray}\label{F3bound}
&&\Xi^{(0,0,12)}_{3;(t_1-s_1,t_2-s_1,s_2-s_1,t_3-s_1)}\vp(x)\nonumber\\
&&\qquad\leq
C\int db\, p_{\io(s_2-s_1)}(x_3,b)\nonumber\\[-8pt]\\[-8pt]
&&\qquad\quad\hspace*{16.8pt}{}\times \int da_1\,da_2\,da_3\,da_4\, p_{\io
(t_1-s_1)}(x_1
,a_1)p_{\io(t_2-s_1)}(x_2,a_2)\nonumber\\
&&\qquad\quad\hspace*{37pt}{} \times
p_{\io(t_3-s_2)}(b,a_3)p_{\io(t_3-s_2)}(b,a_4)\vp
(a_1,a_2,a_3,a_4)\nonumber
\end{eqnarray}
for all $x\in\R^{3\times d}$, $0\leq s_1\leq t_1, t_2\leq s_2\leq t_3$.
Now, from inequality (\ref{F3bound}),
\begin{eqnarray*}
&&\int_{t_2}^{t_3} ds \bigl\langle\Xi^{(0,0,12)}_{3;(t_1,t_2,s,t_3)}\vp
,\mu_0^3\bigr\rangle\\
&&\qquad\leq C\int_{t_2}^{t_3} ds \int\mu_0(dx_3) \int dy\, p_{\io
s}(x_3,y)\\
&&\qquad\quad\hspace*{23.23pt}{}\times
\int dz_3\,dz_4\, p_{\io(t_3-s)}(y,z_3)p_{\io(t_3-s)}(y,z_4)\int dz_1 \,\phi(z_1-z_3) \\
&&\qquad\quad\hspace*{23.23pt}{}\times\int dz_2\,
\phi(z_2-z_4) \int\mu_0(dx_1)p_
{\io t_1}(x_1,z_1)\int\mu_0(dx_2)p_{\io t_2}(x_2,z_2)\\
&&\qquad \leq C(T)\Vert\phi\Vert_1^2.
\end{eqnarray*}
It thus follows that
\[
\int_0^{t_3} dt_2 \int_0^{t_2} dt_1 \int_{t_2}^{t_3} ds \bigl\langle\Xi
^{(0,0,12) }_{3;(t_1,t_2,s,t_3)}\vp,\mu_0^3\bigr\rangle\leq
C(T)\Vert\phi\Vert_1^2.
\]

Again from (\ref{F3bound}), we have that
\[
\sst{i}{j}^3\int_{t_2}^{t_3} ds_2\int_0^{t_1} ds_1 \bigl\langle\Xi
^{(ij,0,0,12)}_{2 ;(s_1,t_1,t_2,s_2,t_3)}\vp, \mu_0^2 \bigr\rangle\leq
C(E_1+E_2),
\]
where for $E_1$ and $E_2$ we have
\begin{eqnarray*}
E_1&=&\sum_{k=1}^2\int_{t_2}^{t_3} ds_2\int_0^{t_1}
ds_1 \int\mu_0(dx_1)\int dy\, p_{\io
s_1}(x_1,y)\int dz\, p_{\io(s_2-s_1)}(y,z)\\
&&\hspace*{12.87pt}{} \times\int dw_{3-k}\,dw_{5-k}\,
p_{\io(t_{3-k}-s_1)}(y,w_{3-k})p_{\io(t_3-s_2)}(z,w_{5-k})\\
&&\hspace*{12.87pt}{}\times \int
dw_{k+2}\,
p_{\io(t_3-s_2)}(z,w_{k+2}) \int dw_k\,
\phi(w_1-w_3)\phi(w_2-w_4) \\
&&\hspace*{12.87pt}{}\times\int\mu_0(dx_2)p_{\io
t_k}(x_2,w_k)\\
&\leq& C\Vert\phi\Vert_1^2\int_{t_2}^{t_3} ds_2
\int_0^{t_1} ds_1 (t_3-s_1)^{-d/2}
\end{eqnarray*}
and
\begin{eqnarray*}
E_2&=&\int_{t_2}^{t_3} ds_2\int_0^{t_1} ds_1 (t_3-s_2)^{-d/2}
\\[-2pt]
&&{}\times\int\mu_0(dx_1)\mu_0(dx_2)\int dw_1\,dw_3\, \phi(w_1-w_3)\\[-2pt]
&&{}\times\int dy \,p_{\io s_1}(x_1,y)p_{\io(t_1-s_1)}(y,w_1) \int
dz \,p_{\io
s_2}(x_2,z)p_{\io(t_3-s_2)}(z,w_3)\\[-2pt]
&&{}\times
\int dw_4 \,p_{\io(t_3-s_2)}(z,w_4)\int dw_2 \,\phi(w_2-w_4)\\[-2pt]
&\leq&
C\Vert\phi\Vert_1^2\int_{t_2}^{t_3} ds_2\int_0^{t_1} ds_1 (t_3-s_2)^{-d/2}.
\end{eqnarray*}
Since $d\leq3$, it follows that
\[
\sst{i}{j}^3\int_0^{t_3} dt_2 \int_0^{t_2} dt_1 \int_{t_2}^{t_3}
ds_2\int_0^{t_1} ds_1 \bigl\langle\Xi^{(ij,0,0,12)}_{2
;(s_1,t_1,t_2,s_2,t_3)}\vp,\mu_0^2\bigr\rangle\leq C(T)\Vert\phi\Vert_1^2.
\]

With one final application of (\ref{F3bound}), we have
\[
\sst{i}{j}^3\int_{t_2}^{t_3} ds_3 \int_0^{t_1} ds_2 \int
_0^{s_2} ds_1
\bigl\langle\Xi^{(12,ij,0,0,12)}_{1;(s_1,s_2,t_1,t_2,s_3,t_3)}\vp,\mu
_0\bigr\rangle\leq
C(F_1+F_2+F_3),
\]
where, for $F_1$ and $F_2$, we have
\begin{eqnarray*}
F_1&=&\int_{t_2}^{t_3} ds_3 \int_0^{t_1} ds_2 \int_0^{s_2} ds_1(t_2-s_1)^{
-d/2}
\int\mu_0(dx) \int dz\, p_{\io s_2}(x,z)\\[-2pt]
&&{}
\times\int dw \,p_{\io(s_3-s_2)}(z,w)\int dv_1\,dv_3 \,p_{\io
(t_1-s_2)}(z,v_1)p_{
\io(t_3-s_3)}(w,
v_3)\phi(v_1-v_3)\\[-2pt]
&&{}\times\int dv_4\, p_{\io(t_3-s_3)}(w,v_4)\int dv_2\, \phi
(v_2-v_4)\\[-2pt]
&\leq& C\Vert\phi\Vert_1\int_{t_2}^{t_3} ds_3 \int_0^{t_1} ds_2
\int_0^{s_2}
ds_1(t_2-s_1)^{-d/2} \int\mu_0(dx) \int dz\, p_{\io s_2}(x,z)\\[-2pt]
&&{}
\times\int dv_1\,dv_3\, p_{\io(t_1-s_2)}(z,v_1)p_{\io(t_3-s_2)}(z,
v_3)\phi(v_1-v_3)
\end{eqnarray*}
and
\begin{eqnarray*}
F_2&=&\int_{t_2}^{t_3} ds_3 \int_0^{t_1} ds_2 \int_0^{s_2}
ds_1(t_2-s_2)^{
-d/2}
\int\mu_0(dx) \int dy\, p_{\io s_1}(x,y)\\[-2pt]
&&{}
\times\int dw\, p_{\io(s_3-s_1)}(y,w)\int dv_1\,dv_3\, p_{\io
(t_1-s_1)}(y,v_1)p_{
\io(t_3-s_3)}(w,
v_3)\phi(v_1-v_3)\\[-2pt]
&&{}\times\int dv_4\, p_{\io(t_3-s_3)}(w,v_4)\int dv_2\,
\phi(v_2-v_4)\\[-2pt]
&\leq& C\Vert\phi\Vert_1\int_{t_2}^{t_3} ds_3 \int_0^{t_1} ds_2
\int_0^{s_2}
ds_1(t_2-s_2)^{-d/2} \int\mu_0(dx) \int dy \,p_{\io s_2}(x,y)\\[-2pt]
&&{}
\times\int dv_1\,dv_3\, p_{\io(t_1-s_2)}(y,v_1)p_{\io(t_3-s_1)}(y,
v_3)\phi(v_1-v_3).
\end{eqnarray*}
Therefore, since $d\leq3$,
\begin{eqnarray*}
&&\sst{i}{j}^3 \int_0^{t_3} dt_2 \int_0^{t_2} dt_1
\int_{t_2}^{
t_3 }
ds_3 \int_0^{t_1} ds_2 \int_0^{s_2} ds_1\bigl\langle\Xi^{(12,ij,0 , 0 ,
12)}_{1;(s_1,s_2,t_1,t_2,s_3,t_3)}\vp, \mu_0\bigr\rangle\\
&&\qquad\leq C(T)\Vert
\phi\Vert_1^2.
\end{eqnarray*}

As a total count of the original fourteen $J_k$, the desired bound has
now been
shown for ten. We continue now with
%
%
\begin{eqnarray}\label{F4bound}
&&\Xi^{(0,12,ij,0)}_{2;(t_1-s_1,s_2-s_1,s_3-s_1,t_2-s_1,t_3-s_1)}
\vp(x)\nonumber\\
&&\qquad \leq C \int db_1\,db_2\,
p_{\io(s_2-s_1)}(x_2,b_1)p_{\io(s_3-s_2)}(b_1,
b_2) \nonumber\\[-8pt]\\[-8pt]
&&\qquad\quad\hspace*{17pt}{}\times\int da_1\,da_2\,da_3\,da_4\, p_{\io
(t_1-s_1)}(x_1,a_1)
p_{\io(t_{(7-i-j)\mn3}-s_2)}(b_1,a_{7-i-j })\nonumber\\
&&\qquad\quad\hspace*{21pt}\hspace*{17pt}{} \times p_{\io
(t_{i+1}-s_3)}(b_2,a_{i+1})p_
{\io(t_3-s_3)}(b_2,a_{j+1})\vp(a)\nonumber
\end{eqnarray}
for any $x\in\R^{2\times d}$, $0\leq s_1\leq t_1\leq s_2\leq s_3\leq
t_2$, and
$i,j=1,2,3$, $i< j$.
Applying inequality (\ref{F4bound}) gives
\begin{eqnarray*}
&&\sst{i}{j}^3\int_{t_1}^{t_2} ds_2 \int_{t_1}^{s_2} ds_1 \bigl\langle\Xi
^{(0,12,ij, 0)}_{2;(t_1,s_1,s_2,t_2,t_3)}\vp,\mu^2_0\bigr\rangle\\
&&\qquad \leq
C\Vert\phi\Vert_1^2\sum_{k=1}^2\int_{t_1}^{t_2} ds_2 \int_{t_1}^{s_2}
ds_1 \int\mu_0(dx_2) \int dy\, p_{\io s_k}(x_2,y)\\
&&\qquad\quad\hspace*{45.2pt}{}\times\int dw_2\,dw_4\, p_{\io(t_2-s_k)}(y,w_2)p_{\io(t_3-s_k)}(y,
w_4)\phi(w_2-w_4)
\\
&&\qquad \leq
C\Vert\phi\Vert_1^2\int_{t_1}^{t_2} ds_2\int_{t_1}^{s_2} ds_1
[(t_3-s_1)^{
-d/2}+(t_3-s_2)^{-d/2}].
\end{eqnarray*}
Therefore, since $d\leq3$, we have
\[
\sst{i}{j}^3 \int_0^{t_3} dt_2 \int_0^{t_2} dt_1 \int_{t_1}^{t_2}
ds_2 \int_{t_1}^{s_2} ds_1 \bigl\langle\Xi^{(0,12,ij,
0)}_{2;(t_1,s_1,s_2,t_2,t_3)}\vp,\mu^2_0\bigr\rangle\leq C(T)\Vert\phi
\Vert_1^2.
\]

With a second and final application of (\ref{F4bound}), it follows that
\[
\sst{i}{j}^3\int_{t_1}^{t_2} ds_3 \int_{t_1}^{s_3} ds_2 \int_0^{t_1}
ds_1\bigl\langle\Xi^{(12,0,12,ij,0)}_{1;(s_1,t_1,s_2,s_3,t_2,t_3)}\vp
,\mu_0\bigr\rangle\leq
C(G_1+G_2+G_3),
\]
where for $G_1$, $G_2$ and $G_3$ we have
\begin{eqnarray*}
G_1&=&\int_{t_1}^{t_2} ds_3 \int_{t_1}^{s_3} ds_2 \int_0^{t_1}
ds_1
(t_3-s_2)^{
-d/2}\int\mu_0(dx)\int dy\, p_{\io s_1}(x,y)\\[-3pt]
&&{} \times
\int dv_3 \int dz\, p_{\io(s_2-s_1)}(y,z)\int dw\, p_{\io
(s_3-s_2)}(z,w)p_{
\io(t_3-s_3)}(w,v_3)\\[-3pt]
&&{}\times\int dv_1\, p_{\io(t_1-s_1)}(y,v_1)\phi(v_1-v_3)\int dv_2\,
p_{\io(s_3-t_1)}(w,v_2)\\[-3pt]
&&{}\times\int dv_4\, \phi(v_2-v_4)\\[-3pt]
&\leq& C
\Vert\phi\Vert_1\int_{t_1}^{t_2} ds_3 \int_{t_1}^{s_3} ds_2 \int
_0^{t_1} ds_1
(t_3-s_2)^{-d/2}\int\mu_0(dx) \int dy\, p_{\io s_1}(x,y)\\[-3pt]
&&{}\times\int dv_1\, p_{\io(t_1-s_1)}(y,v_1)\int dv_3\, p_{\io(t_3-s_1)}(y,
v_3)\phi(v_1-v_3),
\\[-3pt]
G_2&=&\int_{t_1}^{t_2} ds_3 \int_{t_1}^{s_3} ds_2 \int_0^{t_1}
ds_1
(t_3-s_3)^{
-d/2}\int\mu_0(dx)\int dy\, p_{\io s_1}(x,y)\\[-3pt]
&&{}
\times\int dv_1\, p_{\io(t_1-s_1)}(y,v_1)\int dv_3\, \phi
(v_1-v_3)\\[-3pt]
&&{}\times \int dz\,
p_{\io(s_2-s_1)}(y,z)p_{\io(t_3-s_2)}(z,v_3)\\[-3pt]
&&{}\times\int dw\, p_{\io(s_3-s_2)}(z,w)\int dv_2\,
p_{\io(t_2-s_3)}(w,v_2)\int dv_4\, \phi(v_2-v_4)\\[-3pt]
&\leq& C
\Vert\phi\Vert_1\int_{t_1}^{t_2} ds_3 \int_{t_1}^{s_3} ds_2 \int
_0^{t_1} ds_1
(t_3-s_3)^{-d/2}\int\mu_0(dx) \int dy\, p_{\io s_1}(x,y)\\[-3pt]
&&{}\times\int dv_1\, p_{\io(t_1-s_1)}(y,v_1)\int dv_3\, p_{\io(t_3-s_1)}(y,
v_3)\phi(v_1-v_3)
\end{eqnarray*}
and
\begin{eqnarray*}
G_3&=&\int_{t_1}^{t_2} ds_3 \int_{t_1}^{s_3} ds_2 \int_0^{t_1}
ds_1
(t_3-s_3)^{
-d/2}\int\mu_0(dx)\int dy\, p_{\io s_1}(x,y)\\[-3pt]
&&{} \times
\int dv_3 \int dz\, p_{\io(s_2-s_1)}(y,z)\int dw\, p_{\io
(s_3-s_2)}(z,w)p_{
\io(t_3-s_3)}(w,v_3)\\[-3pt]
&&{} \times\int dv_1\, p_{\io(t_1-s_1)}(y,v_1)\phi(v_1-v_3)\int dv_2\,
p_{\io(t_3-s_2)}(z,v_2)\\
&&{}\times\int dv_4\, \phi(v_2-v_4)\\[-3pt]
&\leq&
C\Vert\phi\Vert_1^2 \int_{t_1}^{t_2} ds_3 \int_{t_1}^{s_3}
ds_2 \int_0^ { t_1 }
ds_1 (t_3-s_1)^{-d/2}\\[-3pt]
&&\hspace*{117.1pt}{}\times
\bigl((t_3-s_2)^{-d/2}+(t_3-s_3)^{-d/2}\bigr).
\end{eqnarray*}
Thus
\begin{eqnarray*}
&&G_1+G_2+G_2\\
&&\qquad\leq
C\Vert\phi\Vert_1^2\int_{t_1}^{t_2} ds_3 \int_{t_1}^{s_3}
ds_2 \\
&&\qquad\quad{}\times\int_0^ { t_1 }
ds_1 (t_3-s_1)^{-d/2}\bigl((t_3-s_2)^{-d/2}+(t_3-s_3)^{-d/2}\bigr).
\end{eqnarray*}
And so, since $d\leq3$,
\[
\sst{i}{j}^3 \int_0^{t_3} dt_2 \int_0^{t_2} dt_1
\int_{t_1}^{
t_2 }
ds_3 \int_{t_1}^{s_3} ds_2 \int_0^{t_1}
ds_3\bigl\langle\Xi^{(12, 0 , 12 , ij , 0) } _ { 1;(
s_1,t_1,s_2,s_3,t_2,t_3)}, \mu_0 \bigr\rangle\leq
C(T)\Vert\phi\Vert_1^2.
\]

It thus remains to show the desired bound on two of the fourteen
original~$J_k$.
As in the previous steps, the bounds will result from the following simpler
bound:
%
%
\begin{eqnarray}\label{F5bound}
&&
\Xi^{(0,12,0,12)}_{2;(t_1-s_1,s_2-s_1,t_2-s_1,s_3-s_1,t_3-s_1)}\vp
(x)\nonumber\\
&&\qquad \leq
C\int db_1\, p_{\io(s_2-s_1)}(x_2,b_1) \int db_2\, p_{\io(s_3-s_2)}(b_1,
b_2) \nonumber\\[-8pt]\\[-8pt]
&&\qquad\quad{} \times\int da_1\,da_2\,da_3\,da_4\,
p_{\io(t_1-s_1)}(x_1,a_1)p_{\io(t_2-s_2)}(b_1,a_2)p_{\io
(t_3-s_3)}(b_2,a_3)\nonumber\\
&&\hspace*{21.22pt}\qquad\quad{} \times p_{
\io(t_3-s_3)}(b_2,a_4)\vp(a)\nonumber
\end{eqnarray}
for any $x\in\R^{2\times d}$, $0\leq s_1\leq t_1\leq s_2\leq t_2\leq
s_3\leq
t_3$. Using inequality (\ref{F5bound}), it follows that
\begin{eqnarray*}
&&\int_{t_2}^{t_3} ds_2 \int_{t_1}^{t_2} ds_1 \bigl\langle\Xi
^{(0,12,0,12)}_{ 2;(t_1,s_1,t_2,s_2,t_3)}\vp,\mu_0^2\bigr\rangle\\
&&\qquad \leq
C\Vert\phi\Vert_1\int_{t_2}^{t_3} ds_2 \int_{t_1}^{t_2}
ds_1 \int\mu_0(dx_2) \\
&&\qquad\quad{}\times\int dy\, p_{\io s_1}(x_2,y)
\int dw_2 \,p_{\io(t_2-s_1)} (y,w_2)\\
&&\qquad\quad{}\times\int dw_4\, p_{\io(t_3-s_1)}(y,w_4)\phi(w_2-w_4)\\
&&\qquad \leq
C\Vert\phi\Vert_1^2\int_{t_2}^{t_3} ds_2 \int_{t_1}^{t_2} ds_1
(t_3-s_1)^{-d/2}
.
\end{eqnarray*}
And so, since $d\leq3$,
\[
\int_0^{t_3} dt_2\int_0^{t_2} dt_1 \int_{t_2}^{t_3} ds_1 \int
_{t_1}^{t_2}
ds_2 \bigl\langle\Xi^{(0,12,0,12)}_{ 2;(t_1,s_1,t_2,s_2,t_3)}\vp,\mu
_0^2\bigr\rangle\leq C(T)\Vert\phi\Vert_1^2.
\]

Finally, once again by (\ref{F5bound}),
\begin{eqnarray*}
&&\int_{t_2}^{t_3} ds_3\int_{t_1}^{t_2} ds_2\int_0^{t_1} ds_1
\bigl\langle\Xi^{(12,0 , 12 ,0,12)}_{1;(s_1,t_1,s_2,t_2,s_3,t_3)}\vp,\mu
_0\bigr\rangle\\
&&\qquad \leq C\Vert\phi\Vert_1^2\int_{t_2}^{t_3} ds_3
\int_{t_1}^{t_2} ds_2\int_0^{t_1} ds_1 (t_2-s_1)^{-d/2}(s_3-s_2)^{-d/2}.
\end{eqnarray*}
It thus follows that
\[
\int_0^{t_3} dt_2 \int_0^{t_2} dt_1 \int_{t_2}^{t_3}
ds_3 \int_ { t_1 } ^ { t_2 }
ds_2 \int_0^{t_1} ds_1\bigl\langle\Xi^{(12,0 , 12
,0,12)}_{1;(s_1,t_1,s_2,t_2,s_3,t_3)}\vp,\mu_0\bigr\rangle\leq C(T)\Vert
\phi\Vert_1^2.
\]
Therefore, from the bounds established above for each $J_k$,
$k=1,\ldots,14$, it
follows that
\[
\int_0^{t_3} dt_2 \int_0^{t_2} dt_1\, \E\langle\vp,\mu_{t_1}\mu
_{t_2}\mu_{t_3} ^2\rangle\leq C(T)\Vert\phi\Vert_{L^1}^2.
\]
\upqed\end{pf*}

\section{\texorpdfstring{Proof of Theorem \lowercase{\protect\ref{thmsemimartingale}}}{Proof of Theorem 4.1}}\label{appD}

From Doob's maximal inequality for martingales and Theorem \ref{thmsp} we
have that for $\phi\in C^\infty_K(\R^d)$, $0\leq T <\infty$,
\[
\E\Bigl(\sup_{0\leq t\leq T}\mu_t(\phi)\Bigr)^2\leq2\mu_0(\phi
)^2+8\E
Z_T(\phi)^2+2\E\biggl(\int_0^T ds \,\mu_s(L\phi)\biggr)^2.
\]
For the second term, from Lemma \ref{lmmoment12sp} and H\"{o}lder's inequality,
\begin{eqnarray*}
\E Z_T(\phi)^2&=&\int_0^T ds\, \E\mu_s(\phi^2)+\int_0^T ds\, \E
\mu_s^2(\Lm\phi)\\
& = &\int_0^T ds\, \mu_0(Q_s\phi^2)+
\int_0^T ds\, \mu_0^2(Q^2_s\Lm\phi)\\
&&{} + \int_0^T ds_1 \int
_0^{s_1} ds_2\,
\mu_0(Q_
{s_2}\Phi_{12}Q^2_{s_1-s_2}\Lm\phi)\\
& \leq&
\Vert m\Vert_\infty\int_0^T ds\, \Vert Q_s\phi^2\Vert_{L^1}+\Vert
m\Vert_\infty^2\int
_0^T ds\, \Vert Q^2_s\Lm\phi\Vert_{L^1}\\
&&{} +
\Vert m\Vert_\infty\int_0^T ds_1 \int_0^{s_1} ds_2\, \Vert\Phi
_{12}Q^2_{s_1-s_2}\Lm\phi\Vert_{L^1}\\
& \leq&
C(T)\Vert\phi\Vert_{L^2}^2+C(T)\sum_{i,j=1}^d\Vert\pa_i\phi
\Vert_{L^1}\Vert\pa_j\phi\Vert
_{L^1} \\
&&{}+ C(T)\sum_{i,j=1}^d\Vert\pa_i\phi\Vert_{L^2}\Vert\pa
_j\phi\Vert_{L^2},
\end{eqnarray*}
where in the above $\{S_t\dvtx t\geq0\}$ is the Brownian transition semigroup.
With regards to the third term above,
\begin{eqnarray*}
&&\E\biggl(\int_0^T ds\, \mu_s(L\phi)\biggr)^2\\
&&\qquad\leq
\int_0^T ds_1 \int_0^T ds_2(\E\mu_{s_1}(L\phi)^2\E\mu_{s_2}
(L\phi)^2)^{1/2}\\
&&\qquad \leq T^2\sup_{0\leq s\leq T}\biggl(\mu_0^2\bigl(Q^2_s (L\phi\otimes
L\phi)\bigr)+
\int_0^s dr\, \mu_0\bigl(Q_r\Phi_{12}Q^2_{s-r}(L\phi\otimes L\phi
)\bigr)\biggr)\\
&&\qquad \leq T^2\biggl(C\Vert L\phi\otimes L\phi\Vert_{L^1}+\sup_{0\leq s\leq
T}\int_0^s dr \int dy \bigl(S_{\io(s-r)}L\phi\bigr)(y)^2\biggr)\\
&&\qquad \leq
C(T)\sum_{i,j,p,q=1}^d(\Vert\pa_i\,\pa_j\phi\Vert_{L^1}\Vert\pa
_p\,\pa_q\phi\Vert_{
L^1}+\Vert\pa_i\,\pa_j\phi\Vert_{L^2}\Vert\pa_p\,\pa_q\phi\Vert_{L^2}).
\end{eqnarray*}

Therefore,
%
%
\begin{eqnarray}\label{boundsemimartingale}
&& \E\Bigl(\sup_{0\leq t\leq T}\mu_t(\phi)\Bigr)^2\nonumber\\
&&\qquad \leq
C(T)\Vert\phi\Vert_{L^2}^2+\sum_{i,j=1}^d(\Vert\pa_i\phi\Vert
_{L^1}\Vert\pa_j\phi\Vert_{L^1}+\Vert\pa_i\phi\Vert_{L^2}\Vert
\pa_j\phi\Vert_{L^2})\\
&&\qquad\quad{}+\sum_{i,j,p,q=1}^d(\Vert\pa_i\,\pa_j\phi\Vert_{L^1}\Vert\pa
_p\,\pa_q\phi\Vert_{L^1}
+\Vert\pa_i\,\pa_j\phi\Vert_{L^2}\Vert\pa_p\,\pa_q\phi\Vert_{L^2}
).\nonumber
\end{eqnarray}

If $\phi\in S_d$, from Rudin (\citeyear{rudin}), Theorem 7.10,
there exists a Cauchy
sequence $\{\phi_n\}\subset C^\infty_K(\R^d)$ converging to $\phi$
in $S_d$.
Thus, from (\ref{boundsemimartingale}), Chebyshev's inequality and a
subsequence
argument from the Borel--Cantelli lemma [cf. Theorem 4.2.3 of Chung
(\citeyear{chung})],
there is a subsequence $\{\phi_{n_k}\}$ such that~$\mu_t(\phi
_{n_k})$ converges
uniformly in $t\in[0,T]$ to $\mu_t(\phi)$ with probability one. Therefore,~$\mu_t(\phi)$ is an a.s. continuous semimartingale for $S_d$.

\section{\texorpdfstring{Proof of Lemma \lowercase{\protect\ref{ito}}}{Proof of Lemma 4.2}}\label{appE}

Fix $T\geq0$, and set $\phi\in C_K^\infty(\R^d)$, then from It\^
{o}'s lemma [cf.
Ikeda and Watanabe (\citeyear{ikeda})],
\begin{eqnarray*}
\int_0^Tdt \langle\psi\otimes\phi,\mu_t\mu_T\rangle& =&
\int_0^Tdt\langle\psi\otimes\phi,\mu_t\mu_t\rangle+
\int_0^Tdt\int_0^tds \langle\psi\otimes L\phi,\mu_s\mu_t\rangle
\\
&&{} +\int_0^T dZ_t(\phi)\int_0^t ds\langle\psi,\mu_s\rangle.
\end{eqnarray*}
\begin{lemma}
For any $\phi,\psi\in C_K^\infty(\R^d)$,
\[
\int_0^T dZ_t(\phi)\int_0^t ds\,\mu_s(\psi)=\int_0^T\int_{\R^d}
Z(dt,dx)\int_0^t ds\, \phi(x)\mu_s(\psi).
\]
\end{lemma}
\begin{pf}
Let $0\leq t\leq T$. It follows from (\ref{eqSPDE}) and Corollary
\ref{corbound1} that
\begin{eqnarray*}
&&\E\biggl(\int_0^t dZ_s(\phi) \int_0^s dv\,\mu_v(\psi)
\biggr)^2\\
&&\qquad=\E\int_0^t d
\langle Z(\phi)_\cdot\rangle_s \biggl(\int_0^s dv\,
\mu_v(\psi)\biggr)^2\\
&&\qquad =
\int_0^t ds \int_0^s dv_1 \int_0^s dv_2\, \E\mu_s(\phi^2)\mu
_{v_1}(\psi)\mu_
{v_2}(\psi)\\
&&\qquad\quad{}
+\int_0^t ds \int_0^s dv_1 \int_0^s dv_2\, \E\mu_s(\Lm\phi
)^2\mu_{v_1}
(\psi)\mu_{v_2}(\psi)\\
&&\qquad \leq C(T)(\Vert\phi\Vert_\infty^2\Vert\psi\Vert_\infty^2+
\Vert\Lm\phi\Vert_\infty^2\Vert\psi\Vert_\infty^2
).
\end{eqnarray*}
By assumption on $\Lm$ and since $\phi,\psi\in C^\infty_K(\R^d)$,
$\Vert\Lm\phi\Vert_\infty<\infty$ and$\Vert\Lm\psi\Vert
_\infty<\infty$,
which, since $\Vert\psi\Vert_\infty,\Vert\phi\Vert_\infty
<\infty$, implies by the
definition of
the stochastic integral [cf. Karatzas and Shreve (\citeyear{shreve}),
Chapter 3] that
\[
\int_0^t dZ_s(\phi)\int_0^s dv\, \mu_v(\psi)\in L^2(\P).
\]
In addition, it is clear from Lemma \ref{corbound1} that
$\int_0^s dv\, \mu_v(\psi)\in L^2(\P)$, and thus, again from the
definition of
the stochastic integral,
$ {\int_0^s dv\, \mu_v(\psi)}$ can be approximated in
$L^2(\P)$ by
simple functions of the form
\[
\sum_{i}^n\sum_{A_i} c_{A_i}
1_{A^{(n)}_i}(\omega)1_{(t^{(n)}_i,t^{(n)}_{i+1}]}(s),
\]
where $\bigcup_{i}A^{(n)}_i=\Om$, $A^{(n)}_i\cap A^{(n)}_j=\emp$ if
$i\neq j$,
$\bigcup_{i}(t^{(n)}_i,t^{(n)}_{i+1}]=[0,\infty)$, and
$(t^{(n)}_i,t^{(n)}_{i+1}]\cap(t^{(n)}_k,t^{(n)}_{k+1}]=\emp$ if
$i\neq k$.
It follows that an $L^2$ approximation to
$\int_0^t dZ_s(\phi)\int_0^s dv\, \mu_v(\psi)$ is given by
\begin{eqnarray*}
&&\int_0^t dZ_s(\phi)\sum_{i}^n\sum_{A_i} c_{A_i}
1_{A^{(n)}_i}(\omega)1_{(t^{(n)}_i,t^{(n)}_{i+1}]}(s)\\
&&\qquad = \sum_{i}^n\sum_{A_i} c_{A_i}
1_{A^{(n)}_i}(\omega)1_{(0,t]}(t_i)\bigl(Z_{t_{i+1}}(\phi
)-Z_{t_i}(\phi)\bigr).
\end{eqnarray*}
Clearly $f(s,\phi(x))=\phi(x)\int_0^s dv\, \mu_v(\psi)$ is also
in $L^2(\P)$, and so there exist simple functions of the form
\[
\sum_{i}^n\sum_{A_i} c_{A_i}
1_{A^{(n)}_i}(\omega)1_{(t^{(n)}_i,t^{(n)}_{i+1}]}(s)\phi(x),
\]
converging to $f(s,\phi(x))$ in $L^2(\P)$.
From Walsh's construction of the stochastic integral with respect to a
martingale measure [Walsh (\citeyear{walsh})], an $L^2$ approximation to
$\int_0^t\int Z(ds,dx)\phi(x)\int_0^s dv\, \mu_v(\psi)$ is
then given by
\begin{eqnarray*}
&&\sum_{i}^n\sum_{A_i}\int c_{A_i}
1_{A^{(n)}_i}(\omega)1_{(0,t]}(t_i)\phi(x)(Z_{t_{i+1}}-Z_{t_i})(dx)\\
&&\qquad = \sum_{i}^n\sum_{A_i} c_{A_i}
1_{A^{(n)}_j}(\omega)1_{(0,t]}(t_i)\bigl(Z_{t_{i+1}}(\phi
)-Z_{t_i}(\phi)\bigr).
\end{eqnarray*}
Since any two $L^2$ limits of a sequence must agree, it follows that
\[
\int_0^T dZ_t(\phi)\int_0^t ds\,\mu_s(\psi)=\int_0^T\int_{\R^d}
Z(dt,dx)\int_0^t ds\, \phi(x)\mu_s(\psi).
\]
\upqed\end{pf}

Immediately, we arrive at the corollary:
\begin{coro}
%
%
\begin{eqnarray}\label{coritop1}\qquad
\int_0^T dt \int_0^t ds\, \mu_s(\psi)\mu_t(L\phi)& =&
\int_0^Tdt\, \mu_t(\psi)\mu_T(\phi)
-\int_0^T dt\, \mu_t(\psi)\mu_t(\phi)\nonumber\\[-8pt]\\[-8pt]
&&{} -\int_0^T\int_{\R^d} Z(dt,dx)\int_0^t ds\, \phi(x)\mu
_s(\psi)\nonumber
\end{eqnarray}
for any $\psi,\phi\in C_K^\infty(\R^d)$.
\end{coro}

We can now prove the desired lemma.
\begin{pf}
Assume that $\Psi\in S_{2d}$, then from Lemma
\ref{lmsinglevariableapprox} we can choose $\{\Psi_n;n\in\N\}$ such that
$\Psi_n(x,y)=\sum_{k=1}^n(\psi_k\otimes\phi_k)(x,y)$, for some
\mbox{$\{\psi_k\dvtx k\in\N\}$}, $\{\phi_k\dvtx k\in\N\}\subset C_K^\infty(\R^d)$,
and $\Psi_n$
converges to $\Psi$ in $S_{2d}$ as $n\rightarrow\infty$.
It is clear from (\ref{coritop1}) that
%
%
\begin{eqnarray}\label{iton}\qquad
\int_0^T dt \int_0^t ds \langle L_2\Psi_n,\mu_s\mu_t\rangle& =&
\int_0^Tdt \langle\Psi_n,\mu_t\mu_T\rangle
-\int_0^T dt \langle\Psi_n,\mu_t\mu_t\rangle\nonumber\\[-8pt]\\[-8pt]
&&{}
-\int_0^T\int_{\R^d} Z(dt,dy)\int_0^t ds \langle\Psi_n(\cdot
,y),\mu_s\rangle
.\nonumber
\end{eqnarray}

From Corollary \ref{corbound1},
\begin{eqnarray*}
&&\E\biggl\{\int_0^T dt \langle\Psi_n-\Psi_m,\mu_t\mu
_T\rangle\biggr\}^2\\[-2pt]
&&\qquad=\int_0^T dt \int_0^T ds\, \E\langle(\Psi_n-\Psi_m)\otimes(\Psi
_n-\Psi_m), \mu_t\mu_T\mu_s\mu_T\rangle\\[-2pt]
&&\qquad\leq C(T)\Vert\Psi_n-\Psi_m\Vert_\infty^2
\end{eqnarray*}
and
\begin{eqnarray*}
&&\E\biggl\{\int_0^T dt \langle\Psi_n-\Psi_m,\mu_t\mu
_t\rangle\biggr\}
^2\\[-2pt]
&&\qquad
=\int_0^T dt \int_0^T ds\, \E\langle(\Psi_n-\Psi_m)\otimes(\Psi
_n-\Psi_m), \mu_t\mu_t\mu_s\mu_s\rangle\\[-2pt]
&&\qquad\leq C(T)\Vert\Psi_n-\Psi_m\Vert_\infty^2,
\end{eqnarray*}
since $\Psi_n$ converges in $S_d$ to $\Psi$,
$\lim_{n\rightarrow\infty}\Vert\Psi_n-\Psi\Vert_\infty$
and both of the above two terms are $L^2$ convergent.

For any $t,s\geq0$, since $\mu_\cdot\in C_{M_F(\R^d)}[0,\infty)$, and
$\Psi_n\rightarrow\Psi$ uniformly,
$\langle\Psi_n$, $\mu_s\mu_t\rangle\rightarrow\langle\Psi,\mu
_s\mu_t\rangle$ a.s.
Since the $L^2$ limit
must agree with the a.s. limit,
$L^2 - \lim_{n\rightarrow\infty}\langle\Psi_n,\mu_s\mu_t\rangle
=\langle\Psi,\mu_s\mu_t\rangle$.
Thus,
$L^2 - \lim_{n\rightarrow\infty}\int_0^T dt \langle\Psi_n,\mu
_t\mu_T\rangle=\int_0^T dt
\langle\Psi,\mu_t\mu_T\rangle$, and
$L^2 - \lim_{n\rightarrow\infty}\int_0^T dt \langle\Psi_n,\mu
_t\mu_t\rangle=\int_0^T dt
\langle\Psi,\mu_t\mu_t\rangle$.

Consider next the stochastic integral term and the term involving the
generator~$L$.
From Lemma \ref{corbound1} it follows that
\begin{eqnarray*}
&&\E\biggl\{\int_0^T dt\int_0^t ds \langle L_2\Psi_n-L_2\Psi_m,\mu_s\mu
_t\rangle\biggr\}
^2\\[-2pt]
&&\qquad\leq C(T)\Vert L_2(\Psi_n-\Psi_m)\Vert_\infty^2\\[-2pt]
&&\qquad \leq C(T)\sum_{i,j,p,q=1}^d \Vert\pa_{2_i}\pa_{2_j}(\Psi_n-\Psi
_m)\Vert_\infty
\Vert\pa_{2_p}\pa_{2_q}(\Psi_n-\Psi_m)\Vert_\infty
\end{eqnarray*}
and
\begin{eqnarray*}
&&\E\biggl\{\int_0^T \int_{\R^d}
Z(dt,dy) \int_0^t ds \langle\Psi_n(\cdot,y)-\Psi_m(\cdot,y),\mu
_s\rangle\biggr\}
^2\\[-2pt]
&&\qquad\leq C (H_1 + H_2),
\end{eqnarray*}
where $H_1$ and $H_2$ satisfy
\begin{eqnarray*}
H_1&=&
\int_0^T dt \int_0^t ds_1 \int_0^t ds_2\,\E\langle(\Psi_n-\Psi_m)(x-z)(\Psi_n-\Psi_m)(y-z),\\[-2pt]
&&\hspace*{171pt}{}\mu
_{s_1}(dx)\mu_{s_2} (dy)\mu_t(dz)\rangle\\[-2pt]
&\leq& C
\int_0^T dt \int_0^t ds_2
\int_0^{s_2} ds_1 \,\E\langle(\Psi_n-\Psi_m)(x-z)\cdot(\Psi
_n-\Psi_m)(y-z),\\[-2pt]
&&\hspace*{191.65pt}\mu_{s_1}(dx)\mu_{s_2} (dy)\mu_t(dz)\rangle
\end{eqnarray*}
and
\begin{eqnarray*}
H_2&=&\sum_{i,j=1}^d\int_0^T dt \int_0^t ds_2 \int_0^t ds_1\,\E
\langle\pa_i(\Psi_n-\Psi_m)\otimes\pa_j(\Psi_n-\Psi_m),\mu
_{s_1} \mu_t\mu_{s_2} \mu_t\rangle\\
&\leq&
C\sum_{i,j=1}^d\int_0^T dt \int_0^t ds_2 \int_0^{s_2} ds_1\,\E
\langle\pa_i(\Psi_n-\Psi_m)\,\pa_j(\Psi_n-\Psi_m),\mu_{s_1} \mu
_t\mu_{s_2} \mu_t\rangle.
\end{eqnarray*}
Thus,
\begin{eqnarray*}
H_1+H_2
&\leq& C(T)\Vert\Psi_n-\Psi_m\Vert_\infty^2\\
&&{} + C(T)\sum_{i,j=1}^d
\Vert\pa_{2_i}(\Psi_n-\Psi_m)\Vert_\infty\Vert\pa_{2_j}(\Psi
_n-\Psi_m)\Vert_\infty.
\end{eqnarray*}
Lemma \ref{lmsinglevariableapprox} implies $\Psi_n$ converges in the Schwartz
space $S_{2\times d}$, and thus $\{\pa_{2_i}\Psi_n\dvtx n\in\N\}$ and
$\{\pa_{2_i}\pa_{2_j}\Psi_n\dvtx n\in\N\}$, for all $i,j=1,2,\ldots,d$,
are uniformly
Cauchy sequences.

For any $t,s\geq0$, since $\mu_\cdot\in
C_{M_F(\R^d)}[0,\infty)$, and
$D^\al\Psi_n\rightarrow D^\al\Psi$ uniformly for any multiindex
$\al$, $\langle D^\al\Psi_n,\mu_s\mu_t\rangle\rightarrow\langle
D^\al\Psi,\mu_s\mu_t\rangle$ a.s. Since the
$L^2$ limit must agree with the a.s. limit,
$L^2 - \lim_{n\rightarrow\infty}\langle D^\al\Psi_n,\mu_s\mu
_t\rangle=\langle D^\al\Psi, \mu_s\mu_t\rangle$, and so
$L^2-\lim_{n\rightarrow\infty}\int_0^T dt \int_0^t ds \langle
L_2\Psi_n,\mu_s\mu_t\rangle
=\int_0^T dt \int_0^t ds \langle L_2\Psi,\mu_s\mu_t\rangle$.

Finally,
$
\{\int_0^T \int Z(dt,dy) \int_0^t ds \langle(\Psi_n-\Psi_m)(x,y),
\mu_s(dx)\rangle\}$
is a Cauchy sequence in $L^2$.
Now, for each $y\in\R^d$, and $t\in[0,T]$, $\langle\Psi_n(\cdot
,y),\mu_t\rangle$ is
Cauchy in~$L^2$, and so there exists an a.s convergent subsequence
$\langle\Psi_{n_k}(\cdot,y),\mu_t\rangle$. Since~$\mu_t$ is
almost surely
finite and
$\Psi_n\rightarrow\Psi$ uniformly,
$\langle\Psi_{n_k}(\cdot,y),\mu_t\rangle\rightarrow\langle\Psi
(\cdot,y),\mu_t\rangle$
a.s. as $k\rightarrow\infty$. Furthermore, both $\langle\Psi
_n(\cdot,y),\mu_t\rangle$ and
$\langle\Psi(\cdot,y),\mu_t\rangle$
are uniformly continuous in $y\in\R^d$ and $t\in[0,T]$, and so
\begin{eqnarray*}
&&\lim_{k\rightarrow\infty} \int_0^T \int_{\R^d} Z(dt,dy)
\int_0^t ds\langle\Psi_ { n_k } (\cdot,y),\mu_s\rangle\\
&&\qquad= \int
_0^T \int_{\R^d} Z(dt,dy)
\int_0^t ds
\langle\Psi(\cdot, y) ,\mu_s\rangle\qquad\mbox{a.s.}
\end{eqnarray*}
Since the $L^2$ limit must agree with the a.s. limit,
\begin{eqnarray*}
&&
L^2 - \lim_{n\rightarrow\infty} \int_0^T \int_{\R
^d} Z(dt,
dy)\int_0^t ds\langle\Psi_ {n}(\cdot,y),\mu_s\rangle\\
&&\qquad= \int_0^T
\int_{\R^d} Z(dt,
dy) \int_0^t ds\langle\Psi(\cdot,y),\mu_s\rangle.
\end{eqnarray*}
\upqed\end{pf}

\section{\texorpdfstring{Proof of Theorem \lowercase{\protect\ref{thmexistence}}}{Proof of Theorem 4.3}}\label{appF}

Let $\{G_\eps\}\subset C^\infty(\R^d)$ be any sequence such that
$G_\eps$ and
$\pa_i G_\eps$ converge, respectively, in $L^1$ to $\Gr$ and $\pa
_i\Gr$, and for
$\eps_1,\eps_2>0$, $x\in\R^d$, define
\[
\phi_{\eps_1,\eps_2}(x)=G_{\eps_1}(x)-G_{\eps_2}(x).
\]
Then for the two nonstochastic integral terms, it is clear that
%
%
\begin{eqnarray}\label{GSILT1}
&&\E\biggl[\int_0^T dt\int_0^t ds \langle\phi_{\eps_1,\eps_2},\mu
_s\mu_t\rangle\biggr]^2\nonumber\\
&&\qquad\leq
C \int_0^T dt\, \E\biggl[\int_0^t ds \langle\phi_{\eps_1,\eps_2},\mu
_s\mu_t\rangle\biggr]^2\nonumber\\[-8pt]\\[-8pt]
&&\qquad = C
\int_0^T dt\int_0^t ds_1 \int_0^t ds_2 \,\E\langle\phi_{\eps
_1,\eps_2}\otimes\phi_{\eps_1,\eps_2},\mu_{s_1}\mu_{t} \mu
_{s_1}\mu_{t}\rangle\nonumber\\
&&\qquad\leq
C\int_0^T dt_3\int_0^{t_3} dt_2 \int_0^{t_2} dt_1\,\E\langle\vp
_{\eps_1,\eps_2}, \mu_{t_1}\mu_{t_2}\mu_{t_3}\mu_{t_3}\rangle
\nonumber
\end{eqnarray}
and
%
%
\begin{eqnarray}\label{GSILT2}
&&\E\biggl[\int_0^T dt \langle\phi_{\eps_1,\eps_2},\mu_t\mu_T\rangle
\biggr]^2 \nonumber\\
&&\qquad=
\int_0^T dt_1\int_0^T dt_2\, \E\langle\phi_{\eps_1,\eps_2},\mu
_{t_1}\mu_T\mu_{t_2 }\mu_T\rangle\\
&&\qquad \leq
C\int_0^T dt_2\int_0^{t_2} dt_1\, \E\langle\vp_{\eps_1,\eps
_2},\mu_{t_1}\mu_{t_2} \mu_T\mu_T\rangle,\nonumber
\end{eqnarray}
where $\vp_{\eps_1,\eps_2}(x_1,x_2,x_3,x_4)\triangleq\phi_{\eps
_1,\eps_2}(x_1-x_3)\phi_{\eps_1,\eps_2}(x_2-x_4)$.

For the stochastic integral term, we have
%
%
\begin{eqnarray}\label{GSILT3}
&&\E\biggl[ \int_0^T \int_{\R^d}
Z(dt,dy) \int_0^t ds \langle\phi_{\eps_1,\eps_2}(\cdot-y),\mu
_s\rangle\biggr]^2\nonumber\\
&&\qquad
=\int_0^T dt \,\E\biggl\langle\biggl[\int_0^t ds \langle\phi_{\eps_1,\eps
_2}(\cdot-\cdot\cdot),\mu_s(\cdot) \rangle\biggr]^2,\mu_t\biggr\rangle
\nonumber\\
&&\qquad\quad{} +
\int_0^T dt \biggl\langle\Lm\biggl[\int_0^t ds \langle\phi_{\eps_1,\eps
_2}(\cdot-\cdot\cdot),\mu_s(\cdot) \rangle\biggr],\mu_t\mu_t\biggr\rangle
\\
&&\qquad\leq
C \int_0^T dt_3 \int_0^{t_3} dt_2 \int_0^{t_2}
dt_1
\Biggl[ \E\langle\hat{\vp}_ { \eps_1 , \eps_2},\mu_{t_1}\mu_{t_2}\mu
_{t_3}\rangle\nonumber\\
&&\qquad\quad\hspace*{112.7pt}{} +
\sum_{p,q=1}^d \E\langle\vp^{pq}_{\eps_1,\eps_2},\mu_{t_1}\mu
_{t_2}\mu_{t_3} \mu_ { t_3}\rangle\Biggr],\nonumber
\end{eqnarray}
where
\[
\vp^{pq}_{\eps_1,\eps_2}(x_1,x_2,x_3,x_4)\triangleq
\pa_p\phi_{\eps_1,\eps_2}(x_1-x_3)\,\pa_q\phi_{\eps_1,\eps_2}(x_2-x_4),
\]
for each $p,q=1,2,\ldots,d$, and $(x_1,x_2,x_3,x_4)\in\R^{4\times d}$,
and
\[
\hat{\vp}_{\eps_1,\eps_2}(x_1,x_2,x_3)\triangleq\vp_{\eps_1,\eps
_2}(x_1,x_2,x_3,x_3),
\]
for each $(x_1,x_2,x_3)\in\R^{3\times d}$. Thus, from
(\ref{GSILT1}), (\ref{GSILT2}), (\ref{GSILT3}) and Lemma
\ref{momentbound}, it follows that
\begin{eqnarray*}
\E\biggl[\int_0^T dt\int_0^t ds \langle\phi_{\eps_1,\eps_2},\mu_s\mu
_t\rangle\biggr]^2 &\leq&
C(T)\Vert\phi_{\eps_1,\eps_2}\Vert_1^2,
\\
\E\biggl[\int_0^T dt \langle\phi_{\eps_1,\eps_2},\mu_t\mu_T\rangle
\biggr]^2&\leq& C(T)\Vert\phi_{\eps_1,\eps_2}\Vert_1^2
\end{eqnarray*}
and
\begin{eqnarray*}
&&\E\biggl[ \int_0^T \int_{\R^d}
Z(dt,dy) \int_0^t ds \langle\phi_{\eps_1,\eps_2}(\cdot-y),\mu
_s\rangle\biggr]^2\\
&&\qquad \leq C(T)\Vert\phi_{\eps_1,\eps_2}\Vert_1^2+C(T)\sum
_{p,q=1}^d\Vert\pa_p \phi_{\eps_1,\eps_2}\Vert_1\Vert\pa_q\phi
_{\eps_1,\eps_2}\Vert_1.
\end{eqnarray*}
Since $G_\eps$ and $\pa_i G_\eps$ converge, respectively, to $\Gr$
and $\pa_i\Gr$
in $L^1$, $i=1,\ldots,d$, we have that
$\lim_{\eps_1,\eps_2\rightarrow0}\Vert\phi_{\eps_1,\eps
_2}\Vert_1=0$ and $\lim_{\eps_1,\eps_2\rightarrow
0}\Vert\pa_i\phi_{\eps_1,\eps_2}\Vert_1=0$.

Since the choice of the $\{G_\eps\}$ is arbitrary, it may be assumed that
$G_\eps=\G$ for each $\eps>0$, and the result follows.
\end{appendix}

\section*{Acknowledgment}

Thank you to my advisor Hao Wang, without whom this would
not have been possible.


%

%
\printaddresses

\end{document}